\documentclass[12pt]{article}
\usepackage{amsmath,amsfonts,amssymb,amsthm}
\input amssym.def
\begin{document}
\def \Z{\Bbb Z}
\def \C{\Bbb C}
\def \R{\Bbb R}
\def \Q{\Bbb Q}
\def \N{\Bbb N}
\def \bR{\bf R}
\def \D{{\cal{D}}}
\def \E{{\cal{E}}}
\def \L{\mathcal{L}}
\def \S{{\cal{S}}}
\def \R{{\cal{R}}}
\def \wt{{\rm wt}}
\def \tr{{\rm tr}}
\def \span{{\rm span}}
\def \Res{{\rm Res}}
\def \Der{{\rm Der}}
\def \End{{\rm End}}
\def \Ind {{\rm Ind}}
\def \Irr {{\rm Irr}}
\def \Aut{{\rm Aut}}
\def \Hom{{\rm Hom}}
\def \mod{{\rm mod}}
\def \ann{{\rm Ann}}
\def \ad{{\rm ad}}
\def \rank{{\rm rank}\;}

\def \<{\langle} 
\def \>{\rangle} 
\def \t{\tau }
\def \a{\alpha }
\def \e{\epsilon }
\def \l{\lambda }
\def \La{\Lambda }
\def \g{{\frak{g}}}
\def \h{{\hbar}}
\def \k{{\frak{k}}}
\def \sl{{\frak{sl}}}
\def \gl{{\frak{gl}}}
\def \b{\beta }
\def \om{\omega }
\def \o{\omega }
\def \c{\chi}
\def \ch{\chi}
\def \cg{\chi_g}
\def \ag{\alpha_g}
\def \ah{\alpha_h}
\def \ph{\psi_h}
\def \be{\begin{equation}\label}
\def \ee{\end{equation}}
\def \bex{\begin{example}\label}
\def \eex{\end{example}}
\def \bl{\begin{lem}\label}
\def \el{\end{lem}}
\def \bt{\begin{thm}\label}
\def \et{\end{thm}}
\def \bp{\begin{prop}\label}
\def \ep{\end{prop}}
\def \br{\begin{rem}\label}
\def \er{\end{rem}}
\def \bc{\begin{coro}\label}
\def \ec{\end{coro}}
\def \bd{\begin{de}\label}
\def \ed{\end{de}}

\newcommand{\n}{\:^{\times}_{\times}\:}
\newcommand{\nno}{\nonumber}
\newcommand{\nord}{\mbox{\scriptsize ${\circ\atop\circ}$}}
\newtheorem{thm}{Theorem}[section]
\newtheorem{prop}[thm]{Proposition}
\newtheorem{coro}[thm]{Corollary}
\newtheorem{conj}[thm]{Conjecture}
\newtheorem{example}[thm]{Example}
\newtheorem{lem}[thm]{Lemma}
\newtheorem{rem}[thm]{Remark}
\newtheorem{de}[thm]{Definition}
\newtheorem{hy}[thm]{Hypothesis}
\makeatletter
\@addtoreset{equation}{section}
\def\theequation{\thesection.\arabic{equation}}
\makeatother
\makeatletter

\begin{center}{\Large \bf A smash product construction of
nonlocal vertex algebras}
\end{center}

\begin{center}
{Haisheng Li\footnote{Partially supported by an NSA grant}\\
Department of Mathematical Sciences, Rutgers University, Camden, NJ 08102\\
and\\
Department of Mathematics, Harbin Normal University, Harbin, China}
\end{center}

\begin{abstract}
A notion of vertex bialgebra and a notion of module nonlocal
vertex algebra for a vertex bialgebra are studied and then
a smash product construction of nonlocal vertex algebras is presented.  
For every nonlocal vertex algebra $V$ satisfying a suitable condition, 
a canonical bialgebra $B(V)$ is constructed
such that primitive elements of $B(V)$ are essentially pseudo derivations
and group-like elements are essentially pseudo endomorphisms.
Furthermore, vertex algebras associated with Heisenberg Lie algebras
as well as those associated with nondegenerate even lattices 
are reconstructed through smash products.
\end{abstract}

\section{Introduction}

In the theory of Hopf algebras, there is a notion of $H$-module
algebra for a bialgebra $H$, and furthermore for an $H$-module
algebra $A$, one has the smash product algebra $A\sharp H$.  The
celebrated quantum groups (see [Dr], [J]) are (almost cocommutative)
quasi-triangular Hopf algebras, which generalize Lie algebras and
groups in the sense that the universal enveloping algebras of Lie
algebras and the group algebras of groups are cocommutative Hopf
algebras.  In the theory of Hopf algebras, a very important role has
been played by the Drinfeld's quantum double construction of
quasi-triangular Hopf algebras. It was known (cf. [M]) that if $H$ is
a finite-dimensional cocommutative Hopf algebra, the Drinfeld's
quantum double $D(H)$ is isomorphic to the smash product
$H^{*}\;\sharp \;H$ with respect to a specially defined action of $H$ on
$H^{*}$. Presumably, suitable vertex analogues of the classical smash
product construction and the Drinfeld quantum double construction will
be of importance in the general theory of vertex algebras.

In this paper we formulate and study a notion of vertex
bialgebra and a notion of $B$-module nonlocal vertex algebra for a
vertex bialgebra $B$ and we then give a smash product
construction. The vertex operator algebras associated with Heisenberg
Lie algebras and the vertex algebras associated with nondegenerate
even lattices are studied again in terms of the smash product
construction.  We shall study Drinfeld double construction in a
sequel.

The notion of vertex bialgebra, studied in this paper, is rather
simple, where a vertex bialgebra is just a nonlocal vertex algebra
equipped with a (standard) coalgebra structure such that the
comultiplication map and the counit map are homomorphisms of nonlocal
vertex algebras.  This notion is analogous to and generalizes the
classical notion of bialgebra as the notion of nonlocal vertex algebra
is analogous to and generalizes that of associative algebra.  
Two families of examples of vertex bialgebras are
given in this paper. The first family is associated to standard bialgebras
equipped with a derivation $\partial$ such that $\Delta\circ \partial
=(\partial\otimes 1+1\otimes\partial) \Delta$ and $\varepsilon \circ
\partial =0$, and the second family of examples are based on vertex
algebras associated with vertex Lie algebras ([DLM], [FB], [P]),
namely conformal Lie algebras (see [K]).

Note that a notion of vertex operator coalgebra, which is dual to the
notion of vertex operator algebra, has been formulated and studied by
Hubbard in \cite{hub1} and \cite{hub2}.  We here just use the
classical notion of coalgebra, bypassing the general notion of vertex
operator coalgebra.  Certainly, more general notions of vertex
bialgebra exist.

For a vertex bialgebra $B$, a $B$-module nonlocal vertex algebra is
a nonlocal vertex algebra $V$ equipped with a $B$-module structure
such that the following conditions hold for $b\in B,\; v\in V$:
\begin{eqnarray*}
& &Y(b,x)v\in V\otimes \C((x)),\\
& &Y(b,x){\bf 1}=\varepsilon (b){\bf 1},\\
& &Y(b,x_{1})Y(v,x)=\sum Y(Y(b_{(1)},x_{1}-x)v,x)Y(b_{(2)},x_{1}),
\end{eqnarray*}
where $\Delta(b)=\sum b_{(1)}\otimes b_{(2)}$.  Furthermore, for a
given vertex bialgebra $H$ and an $H$-module nonlocal vertex algebra
$V$, the smash product nonlocal vertex algebra $V\;\sharp\; H$ has the
underlying vector space $V\otimes H$ with the vertex operator map
$Y^{\sharp}$ defined by
$$Y^{\sharp}(u\otimes g,x)(v\otimes h)=\sum Y(u,x)Y(g_{(1)},x)v\otimes
Y(g_{(2)},x)h$$ for $u,v\in V,\; g,h\in H$. 
These notions (or construction) are analogous to and generalize 
their classical counterparts. For this study,
a key role is played by the notions of pseudo derivation and
pseudo endomorphism as we mention next.

In the work \cite{ek} on formal deformations 
of vertex operator algebras, Etingof and Kazhdan introduced 
and studied a notion of pseudo derivation, where
a pseudo derivation of a vertex algebra $V$ is an element $\psi(x)$ of the
space $\Hom_{\C}(V,V\otimes \C((x)))$ such that
\begin{eqnarray}
& &\hspace{1cm}  \ \ \ [L(-1),\psi(x)]=-\frac{d}{dx}\psi(x),
\label{epseudo-der-intro-d}\\
& &\psi(x_{1})Y(v,x_{2})=Y(v,x_{2})\psi(x_{1})+Y(\psi(x_{1}+x_{2})v,x_{2})
\;\;\;\mbox{ for }v\in V.\label{epseudo-der-intro}
\end{eqnarray}
A little bit earlier in \cite{li-sc}, in an attempt to interpret
the physics superselection theory in the context of vertex operator algebras,
we studied objects $\Phi(x)\in \Hom (V,V[x,x^{-1}])$ with
the following properties:
\begin{eqnarray}
& &\Phi(x){\bf 1}={\bf 1},\label{e1pseudo-end-intro}\\
& &[L(-1),\Phi(x)]=-\frac{d}{dx}\Phi(x),\\
& &\Phi(x_{1})Y(v,x)=Y(\Phi(x_{1}+x)v,x)\Phi(x_{1})
\;\;\;\mbox{ for all }v\in V.\label{e3pseudo-end-intro}
\end{eqnarray}
In view of the notion of pseudo derivation, such objects $\Phi(x)$
are nothing but ``pseudo endomorphisms'' and 
this observation was expanded in \cite{li-pseudo}.
What was proved in \cite{li-sc} (cf. \cite{li-twisted}) 
that for such a $\Phi(x)$
and for any $V$-module $(W,Y_{W})$, the
pair $(W,Y_{W}(\Phi(x)\cdot,x))$ carries the structure 
of a $V$-module. If $\Phi(x)$ is invertible, 
it was proved that $(V,Y(\Phi(x)\cdot,x))$ is 
a simple current in a certain sense.
In that study,  $\Phi(x)$ was indeed 
treated as a kind of endomorphism. 

Recall that for a Hopf algebra $H$ with an
$H$-module algebra $A$, primitive elements of $H$ act on $A$ by
derivations while group-like elements of $H$ act on $A$ as
automorphisms. For our particular notions of
vertex bialgebra and module vertex algebra, the roles of
derivations and automorphisms are played by pseudo derivations 
and pseudo automorphisms.
What is more, for any nonlocal vertex algebra $V$ 
satisfying a mild condition,
we construct a canonical bialgebra $B(V)$ 
whose primitive elements are essentially
pseudo derivations and whose group-like elements are essentially pseudo
automorphisms. 
By definition, $B(V)$ is the (unique) maximal subspace of
$\Hom_{\C}(V,V\otimes \C((x)))$, satisfying a certain property.
Note that $\Hom_{\C}(V,V\otimes \C((x)))$,
which can be identified with $\End_{\C((x))}(V\otimes
\C((x)))$, is an associative algebra over $\C((x))$.  
It is proved that for a general $V$, 
$B(V)$ is an associative subalgebra with the formal differential
operator $d/dx$ as a $\C$-linear derivation.  If $V$ is a simple
vertex operator algebra in the sense of \cite{flm}, or if $V$ is a
vertex algebra with a certain basis of P-B-W type, it is proved that
$B(V)$ is furthermore a bialgebra.

Just as a simple associative algebra cannot be made into a bialgebra,
a simple vertex algebra cannot be made into a vertex
bialgebra.  On the other hand, as Borcherds pointed out in
\cite{b-qva}, vertex algebras often have a natural cocommutative
bialgebra structure, and furthermore, the canonical operator $L(-1)$ of
the vertex algebras acts as a derivation.  Then we have a natural
vertex bialgebra.  For example, for the vertex algebra $V_{L}$
associated with a nondegenerate even lattice $L$, the underlying
vector space is the tensor product of a symmetric algebra with the
group algebra of $L$, which is naturally a commutative and
cocommutative Hopf algebra.  This vertex bialgebra is denoted by
$B_{L}$.  Let $\epsilon$
be a $2$-cocycle of $L$, which was used in the construction of the
vertex algebra $V_{L}$ (see \cite{flm}) and let $B_{L,\epsilon}$ 
be the tensor product of
the symmetric algebra with the $\epsilon$-twisted group
algebra of $L$. Naturally, 
$B_{L,\epsilon}$ is a (non-commutative) associative 
algebra with $L(-1)$ as a derivation. 
Then $B_{L,\epsilon}$ becomes a nonlocal vertex algebra.
It is proved that 
there is a natural $B_{L}$-module nonlocal vertex algebra structure
on $B_{L,\epsilon}$ such that the smash product nonlocal
vertex algebra $B_{L,\epsilon}\sharp B_{L}$ contains
$V_{L}$ a subalgebra.
We also obtain similar results for
the vertex operator algebras associated with
Heisenberg Lie algebras, where the vertex operator algebras
are identified with the diagonal subalgebras of the smash products
of type $B\sharp B$ with the comultiplication map as the vertex algebra 
isomorphism. 

This paper is organized as follows: In Section 2, we study pseudo
derivations and pseudo endomorphisms for nonlocal vertex algebras. In
Section 3, we construct a canonical bialgebra
for a nondegenerate nonlocal vertex algebra.
In Section 4, we define the notion of vertex bialgebra and the notion of
module nonlocal vertex algebra for a vertex bialgebra, and we give the
smash product construction.    In Section 5, we use the
smash product construction to construct the vertex algebras associated
with Heisenberg Lie algebras and the vertex algebras associated with
nondegenerate even lattices.

\section{Pseudo derivations and pseudo endomorphisms}
In this section we revisit pseudo-derivations and pseudo-endomorphisms 
for nonlocal vertex algebras, extending some of the results
obtained in \cite{li-pseudo}.

First, we recall the notion of nonlocal 
vertex algebra (\cite{b-Gva}, \cite{bk}, \cite{li-g1}, \cite{li-qva2}).
A {\em nonlocal vertex algebra} is
a vector space $V$ equipped with a linear map
\begin{eqnarray}
Y(\cdot,x): & &V\rightarrow \Hom (V,V((x)))
\subset (\End\; V)[[x,x^{-1}]]\nonumber\\
& &v\mapsto Y(v,x)=\sum_{n\in \Z}v_{n}x^{-n-1}
\;\;\;(\mbox{where }v_{n}\in \End \;V)
\end{eqnarray}
and equipped with a distinguished vector ${\bf 1}$ of $V$, called the
{\em vacuum vector}, such that
\begin{eqnarray}
& & Y({\bf 1},x)=1,\\
& &Y(v,x){\bf 1}\in V[[x]]\;\;\;\mbox{ and }\;\;
\lim_{x\rightarrow 0}Y(v,x){\bf 1}=v\;\;\;\mbox{ for }v\in V
\end{eqnarray}
and such that for $u,v,w\in V$, there exists a nonnegative integer $l$ 
such that
\begin{eqnarray}
(x_{0}+x_{2})^{l}Y(u,x_{0}+x_{2})Y(v,x_{2})w=
(x_{0}+x_{2})^{l}Y(Y(u,x_{0})v,x_{2})w.
\end{eqnarray}

With this notion, an (ordinary) vertex algebra $V$ 
is simply a nonlocal vertex algebra
that satisfies the weak commutativity in the sense that for any $u,v\in V$, 
there exists a nonnegative integer $k$ such that
\begin{eqnarray}
(x_{1}-x_{2})^{k}Y(u,x_{1})Y(v,x_{2})=(x_{1}-x_{2})^{k}Y(v,x_{2})Y(u,x_{1}).
\end{eqnarray}
For a nonlocal vertex algebra $V$,
just as with an ordinary vertex algebra, we define a linear operator 
$\D$ by $\D v=v_{-2}{\bf 1}$ for $v\in V$. Then we have
\begin{eqnarray}
[\D,Y(v,x)]=Y(\D v,x)=\frac{d}{dx}Y(v,x)\;\;\;\mbox{ for }v\in V.
\end{eqnarray}
If $V$ is a vertex operator algebra in the sense 
of \cite{flm} and \cite{fhl}, then $\D =L(-1)$. 

\bex{eborcherds-example}
{\em Let $A$ be any (unital) associative algebra equipped 
with a derivation $\partial$.
Then (see [B1,2]) we have a nonlocal vertex algebra structure on $A$
with $1$ as the vacuum vector and with
\begin{eqnarray}
Y(a,x)b=(e^{x\partial}a)b\;\;\;\mbox{ for }a,b\in A.
\end{eqnarray}
In this case $\D=\partial$. We denote this nonlocal vertex algebra by $(A,\partial)$.
Clearly, $(A,\partial)$ is an ordinary vertex algebra if and only if
$A$ is commutative. On the other hand, let $V$ be a nonlocal vertex algebra 
such that $Y(u,x)v\in V[[x]]$ for all $u,v\in V$. Then
one can prove that $V$ equipped with the multiplicative operation
defined by $u\cdot v=u_{-1}v$ for $u,v\in V$ is
an associative algebra with ${\bf 1}$ as the identity 
and with $\D$ as a derivation.
Furthermore, $Y(u,x)v=(e^{x\D}u)\cdot v$ for $u,v\in V$.}
\eex

\br{rclassical-homo}
{\em  Let $A$ and $B$ be (unital) associative algebras
equipped with derivations $\partial_{A}$ and $\partial_{B}$,
respectively.  Then a linear map $f$ from $A$ to $B$ is a homomorphism
of nonlocal vertex algebras from $(A,\partial_{A})$ to $(B,\partial_{B})$
if and only if $f$ is a homomorphism of
algebras such that $f\partial_{A}=\partial_{B}f$.  }
\er

\br{rclassical-product}
{\em  Let $A$ and $B$ be (unital) associative algebras
equipped with derivations $\partial_{A}$ and $\partial_{B}$,
respectively.  Clearly, $\partial_{A}\otimes 1+1\otimes \partial_{B}$ is 
a derivation of the tensor product associative algebra $A\otimes B$.
Then we have }
$$(A,\partial_{A})\otimes (B,\partial_{B})
=(A\otimes B,\partial_{A}\otimes 1+1\otimes \partial_{B}).$$
\er

Let $V$ and $K$ be nonlocal vertex algebras.
A {\em homomorphism} of nonlocal vertex algebras from  $V$ to $K$ 
is a linear map $\psi$ from $V$ to $K$ such that
\begin{eqnarray}
& &\psi({\bf 1})={\bf 1},\\
& &\psi Y(u,x)v=Y(\psi (u),x)\psi(v)\;\;\;\mbox{ for }u,v\in V.
\end{eqnarray}
It follows that $\psi\D_{V}=\D_{K}\psi$, where $\D_{V}$ and $\D_{K}$ denote 
the $\D$-operators of $V$ and $K$, respectively. 
A homomorphism of nonlocal vertex algebras from $V$ to itself is
called an {\em endomorphism} of $V$.

Let $V$ be a nonlocal vertex algebra. A {\em derivation} of $V$ (see [B1])
is a linear endomorphism $f$ of $V$ such that
\begin{eqnarray}
f Y(v,x)=Y(f(v),x)+Y(v,x)f\;\;\;\mbox{ for }v\in V.
\end{eqnarray}
All the derivations of $V$ form a Lie subalgebra $\Der (V)$ 
of the general linear Lie algebra $\gl (V)$. 
If $V$ is an ordinary vertex algebra, for any $v\in V$, 
$v_{0}$ is a derivation of $V$ and 
$\{ v_{0}\;|\; v\in V\}$ is an ideal of $\Der (V)$.

The following is a simple analogue of a classical notion:

\bd{dsigma-tau-derivation}
{\em Let $V$ and $K$ be nonlocal vertex algebras 
and let $\sigma$ and $\tau$ be
homomorphisms of nonlocal vertex algebras from $V$ to $K$.
A {\em $(\sigma,\tau)$-derivation} from $V$ to $K$ is a linear map 
$\psi$ from $V$ to $K$ such that for $u,v\in V$,
\begin{eqnarray}\label{esigma-phi-der}
\psi(Y(u,x)v)=Y(\sigma(u),x)\psi(v)+Y(\psi(u),x)\tau(v).
\end{eqnarray}
In the case that $V$ is a subalgebra of $K$
with $\sigma$ and $\tau$ being the identity homomorphism,
we simply call a $(\sigma,\tau)$-derivation 
{\em a derivation from $V$ to $K$.} All the derivations from $V$ to $K$
form a subspace which we denote by $\Der (V,K)$.}
\ed

For any derivation $f$ from $V$ to $K$, from (\ref{esigma-phi-der}) we have
\begin{eqnarray}
f({\bf 1})=0,\ \ \ \ \D_{K} f=f\D_{V}.
\end{eqnarray}

Now, let $K$ be the tensor product nonlocal vertex algebra 
$V\otimes (A,\partial)$, where $(A,\partial)$ is 
the nonlocal vertex algebra associated with
an associative algebra $A$ equipped with a derivation $\partial$
(recall Example \ref{eborcherds-example}).
Let $\sigma=\tau$ be the natural embedding of $V$ into $K$.
We are interested in derivations from $V$ to $V\otimes (A,\partial)$.
Consider $V\otimes A$ as an $A$-module in the obvious way. We have
\begin{eqnarray}
\Hom_{\C}(V,V\otimes A)=\End_{A}(V\otimes A)
\end{eqnarray}
as $A$-modules and as $\C$-vector spaces. 
We may and we should consider any linear map from $V$ to $V\otimes A$ 
or from $V$ to $V$ as an $A$-linear endomorphism of $V\otimes A$.
The vertex operator map $Y$ of $V$ is also considered 
as an $A$-linear map from $V\otimes A$ to $(\End (V\otimes A))[[x,x^{-1}]]$, that is,
\begin{eqnarray}
Y(v\otimes a,x)=Y(v,x)\otimes a\;\;\;\mbox{ for }v\in V,\; a\in A.
\end{eqnarray}
On the other hand, we have the vertex operator map which we denote by $Y_{ten}$ for
the tensor product nonlocal vertex algebra $V\otimes (A,\partial)$.

\bl{lgeneral-A}
Let $V$ be a nonlocal vertex algebra and let $A$ be an associative algebra
equipped with a derivation $\partial$. A 
linear map $\psi$ from $V$ to $V\otimes A$
is a derivation from $V$ to the tensor product nonlocal 
vertex algebra $V\otimes (A,\partial)$ if and only if
\begin{eqnarray}
\psi Y(v,x)=Y(v,x)\psi+Y(e^{x(1\otimes \partial)}\psi(v),x)
\;\;\;\mbox{ for }v\in V.
\end{eqnarray}
\el

\begin{proof} For $u\in V,\; a\in A$, we have
$$Y_{ten}(u\otimes a,x)=Y(u,x)\otimes Y(a,x)=
Y(u,x)\otimes (e^{x\partial }a)=Y(e^{x(1\otimes \partial)}(u\otimes a),x).$$
Then
\begin{eqnarray}
Y_{ten}(\psi(v),x)=Y(e^{x(1\otimes \partial)}\psi(v),x)
\;\;\;\mbox{ for }v\in V.
\end{eqnarray}
For $v\in V$, we see that 
$\psi Y(v,x)=Y(v,x)\psi+Y_{ten}(\psi(v),x)$ exactly amounts to
$$\psi Y(v,x)=Y(v,x)\psi+Y(e^{x(1\otimes \partial)}\psi(v),x).$$
This proves the assertion. 
\end{proof}

Now we consider special cases with $A=\C((x))$ or $\C((x))[\log x]$ and with
$\partial=\pm d/dx$.

\bd{dderivation}
{\em Let $V$ be a nonlocal vertex algebra. We define ${\rm PDer}^{\pm}(V)$ to be
the subspaces of $\Hom_{\C}(V,V\otimes \C((x)))$, consisting of elements $\psi(x)$
satisfying the condition
\begin{eqnarray}
[\psi(x_{1}), Y(v,x_{2})]=Y(\psi(x_{1}\pm x_{2})v,x_{2})
\;\;\;\mbox{ for all }v\in V.
\end{eqnarray}}
\ed
 
\bp{ppseudo-minus}
Let $V$ be a nonlocal vertex algebra. We have
\begin{eqnarray}
{\rm PDer}^{\pm}(V)=\Der (V, V\otimes (\C((x)),\pm d/dx)).
\end{eqnarray}
Furthermore, for $\psi(x)\in {\rm PDer}^{\pm }(V)$, we have
\begin{eqnarray}
\psi(x){\bf 1}=0,\ \ \ \ \ \ [\D,\psi(x)]=\mp \frac{d}{dx}\psi(x).
\end{eqnarray}
\ep

\begin{proof} For $f(x)\in \C((x))$, we have
$e^{\pm x_{2}\frac{d}{dx}}f(x)=f(x\pm x_{2})$ (by the Taylor theorem). Then for $v\in V$,
$$Y(e^{\pm x_{2}(1\otimes \frac{d}{dx})}\psi(x)v,x_{2})=Y(\psi(x\pm x_{2})v,x_{2}).$$
Now the first assertion follows from Lemma \ref{lgeneral-A}.
For $\psi(x)\in {\rm PDer}^{\pm }(V)$, since $\psi(x)$ is a derivation from $V$ to
$V\otimes (\C((x)),\pm \frac{d}{dx})$, we have
$\psi(x){\bf 1}=0$ and
$$\psi(x) \D_{V}
=\left(\D\otimes 1\pm 1\otimes \frac{d}{dx}\right)\psi(x),$$
noticing that the $\D$-operator of $V\otimes (\C((x)),\pm \frac{d}{dx})$ is
$\D\otimes 1\pm 1\otimes \frac{d}{dx}$. Then we have
$$[\D,\psi(x)]=\mp \frac{d}{dx} \psi(x)\;\;\mbox{ on }V,$$
completing the proof.
\end{proof}

\br{rek-pseudo}
{\em Recall from \cite{ek} that a {\em pseudo-derivation} of 
an ordinary vertex algebra $V$ is 
a linear map $a(x)$ from $V$ to $V\otimes \C((x))$ such that
\begin{eqnarray}
& &[\D,a(x)]=-\frac{d}{dx}a(x),\\
& &[a(x_{1}), Y(v,x_{2})]=Y(a(x_{1}+x_{2})v,x_{2})
\;\;\;\mbox{ for all }v\in V.
\end{eqnarray}
The space of pseudo-derivations of $V$ was denoted by ${\rm PDer}(V)$.
By Proposition \ref{ppseudo-minus}, a
pseudo-derivation of $V$ in the sense of Etingof and Kazhdan is
exactly a derivation from $V$ to $V\otimes (\C((x)),d/dx)$.}
\er

\br{rlie-algebra}
{\em Let $V$ be any nonlocal vertex algebra. Then
${\rm PDer}^{\pm }(V)$ are subalgebras
of the linear Lie algebra $\gl(V\otimes \C((x)))$ and
the formal differential operator $\frac{d}{dx}$ 
is a $\C$-linear derivation of ${\rm PDer}^{\pm}(V)$. 
It is straightforward to see that
the linear endomorphism of $\Hom (V,V\otimes \C((x)))$, sending
$a(x)$ to $a(-x)$, is a Lie algebra isomorphism between 
${\rm PDer}^{\pm }(V)$. }
\er
 
We have the following result (cf. \cite{ek}):

\bp{pinner-derivation}
Let $V$ be an ordinary vertex algebra. For any 
$u\in V,\; f(x)\in \C((x))[\log x]$, set
$$\Phi^{\pm}(u,f)
=\Res_{x_{1}}e^{\pm x_{1}\frac{\partial}{\partial x}}f(x)Y(u,x_{1})
=\sum_{n\ge 0}\frac{(\pm 1)^{n}}{n!}f^{(n)}(x)u_{n}.$$
Then 
\begin{eqnarray}
\Phi^{\pm}(u,f)\in \Der(V, (\C((x))[\log x],\pm d/dx)).
\end{eqnarray}
Furthermore, if $f(x)\in \C((x))$,
then $\Phi^{\pm}(u,f)\in {\rm PDer}^{\pm}(V)$.
In particular, for $u\in V$,
\begin{eqnarray}
Y(u,x)^{-}=\sum_{n\ge 0}u_{n}x^{-n-1}
=\Phi^{-}(u,x^{-1})\in {\rm PDer}^{-}(V).
\end{eqnarray}
\ep

\begin{proof} For any $v\in V,\; f\in \C((x))[\log x]$, considering
$v\otimes f$ as an element of the vertex algebra 
$V\otimes (\C((x))[\log x],\pm d/dx)$, we have
$$(u\otimes f)_{0}=\Res_{x_{1}}Y_{ten}(u\otimes f,x_{1})
=\Res_{x_{1}}(e^{\pm x_{1}\frac{\partial}{\partial x}}f(x))Y(u,x_{1})=\Phi^{\pm}(u,f).$$
Since $V\otimes (\C((x))[\log x],\pm d/dx)$
are ordinary vertex algebras,
$(u\otimes f)_{0}$ is a derivation of
$V\otimes (\C((x))[\log x],\pm d/dx)$ and hence a derivation from $V$ to
$V\otimes (\C((x))[\log x],\pm d/dx)$. 
If $f(x)\in \C((x))$, by Proposition \ref{ppseudo-minus},
$\Phi^{\pm}(u,f)\in {\rm PDer}^{\mp}(V)$.
\end{proof}

The elements $\Phi^{\pm}(u,f)$ in Proposition \ref{pinner-derivation}
are somewhat inner pseudo derivations.  If $V$ is a nonlocal vertex
algebra associated to a differential associative algebra, such pseudo
derivations are all zero. In this case, we have the following result:

\bp{pclassical-3}
Let $A$ be an associative algebra equipped with a derivation $\partial$.
Suppose that $d(x)\in (\Der A)[[x,x^{-1}]][\log x]$ 
such that 
\begin{eqnarray}
& &d(x)a\in A\otimes \C((x))[\log x]\;\;\;\mbox{ for }a\in A,\\
& &[\partial, d(x)]=\pm \frac{d}{dx} d(x).
\end{eqnarray}
Then $d(x)\in {\rm Der}(A, A\otimes (\C((x))[\log x],\mp \frac{d}{dx}))$.
Furthermore, if $d(x)\in (\Der A)[[x,x^{-1}]]$, we have
$d(x)\in {\rm PDer}^{\mp}(A,\partial)$.
\ep

\begin{proof} {}From the assumption we have
$$e^{x_{1}\partial}d(x)e^{-x_{1}\partial}=e^{\pm x_{1}\frac{d}{dx}}d(x).$$
For $a,b\in A$, we have
\begin{eqnarray*}
d(x)Y(a,x_{2})b&=&d(x)\left((e^{x_{2}\partial}a)b\right)\\
&=&(d(x)e^{x_{2}\partial}a)b+(e^{x_{2}\partial}a)d(x)b\\
&=&\left(e^{x_{2}\partial}e^{\mp x_{2}\frac{d}{dx}}d(x)a\right)b+(e^{x_{2}\partial}a)d(x)b\\
&=&Y(e^{\mp x_{2}\frac{d}{dx}}d(x)a,x_{2})b+Y(a,x_{2})d(x)b.
\end{eqnarray*}
This proves $d(x)\in {\rm Der}(A, A\otimes (\C((x))[\log x],\mp
\frac{d}{dx}))$. The last assertion
follows from Proposition \ref{ppseudo-minus}.
\end{proof}

Next we study endomorphisms.
We have the following straightforward analogue of Lemma \ref{lgeneral-A}:

\bl{lgeneral-pseudo-end}
Let $V$ be a nonlocal vertex algebra and let $A$ be an associative algebra 
equipped with a derivation $\partial$.
A linear map $\phi$ from $V$ to $V\otimes A$ is a 
homomorphism of nonlocal vertex algebras from $V$ to the tensor product nonlocal 
vertex algebra $V\otimes (A,\partial)$ if and only if 
\begin{eqnarray}
& &\phi ({\bf 1})={\bf 1},\\
& &\phi Y(v,x)=Y(e^{x(1\otimes \partial)}\phi(v),x_{2})\phi
\;\;\;\mbox{ for }v\in V. 
\end{eqnarray}
\el

\bd{dpseudo-end}
{\em Let $V$ be a nonlocal vertex algebra. 
We define ${\rm PEnd}^{\pm}(V)$ to be the subspaces of
$\Hom _{\C}(V,V\otimes \C((x)))$, consisting of 
elements $\phi(x)$ satisfying the conditions
\begin{eqnarray}
& &\phi(x){\bf 1}={\bf 1},\\
& &\phi(x_{1})Y(v,x_{2})=Y(\phi(x_{1}\pm x_{2})v,x_{2})\phi(x_{2})
\;\;\;\mbox{ for }v\in V.
\end{eqnarray}}
\ed

The following is an analogue of Proposition \ref{ppseudo-minus}:

\bp{ppseudo-end-diff}
Let $V$ be a nonlocal vertex algebra and let $\phi(x)\in \Hom_{\C}(V,V\otimes \C((x)))$.
Then $\phi(x)\in {\rm PEnd}^{\pm}(V)$ if and only if $\phi(x)$
is a homomorphism of nonlocal vertex algebras from $V$ to
$V\otimes (\C((x)),\pm d/dx)$. Furthermore, 
\begin{eqnarray}
[\D,\phi(x)]=\mp \frac{d}{dx}\phi(x)\;\;\;\mbox{ for }\phi(x)\in {\rm PEnd}^{\pm} (V).
\end{eqnarray}
\ep

The following is an explicit construction of pseudo endomorphisms
(cf. \cite{li-sc}):

\bp{ppseudo-end}
Let $V$ be an ordinary vertex algebra and let $h\in V$ be such that
\begin{eqnarray}
[h(m),h(n)]=0\;\;\;\mbox{ for }m,n\ge 0,
\end{eqnarray}
where $h(k)=h_{k}$ for $k\in \Z$,
and such that $h(0)$ acts semisimply on $V$ with integral eigenvalues.
Then $E^{+}(-h,x)x^{h(0)}\in {\rm PEnd}^{-}(V)$, where
$$E^{+}(-h,x)=\exp \left(\sum_{n\ge 1}\frac{-h(n)}{n}x^{-n}\right).$$
\ep

\begin{proof} Set 
$$\phi(x)=h(0) \log x -\sum_{n\ge 1} \frac{h(n)}{n}x^{-n}
\in \Hom (V,V\otimes \C((x))[\log x]).$$
Then
$$\phi(x)= \sum_{n\ge 0}\frac{(-1)^{n}}{n!}h(n)\left(\frac{d}{dx}\right)^{n}(\log x)
=\Phi^{-} (h,\log x).$$
By Proposition \ref{pinner-derivation}, $\phi(x)$
is a derivation from $V$ to $V\otimes (\C((x))[\log x],-d/dx)$.
Then $E^{+}(-h,x)x^{h(0)}=e^{\phi(x)}$ is an endomorphism from $V$ to
$V\otimes \C((x))\subset V\otimes \C((x))[\log x]$.
By Proposition \ref{ppseudo-end-diff} we have  $E^{+}(-h,x)x^{h(0)}\in {\rm PEnd}^{-}(V)$.
\end{proof}

\section{Differential bialgebra $B(V)$}

In this section we associate a canonical associative algebra $B(V)$ to
each nonlocal vertex algebra $V$ and we prove that $B(V)$ has a
natural bialgebra structure if $V$ satisfies a certain
condition. 

First we formulate the following notion:

\bd{ddelta-closed-subset}
{\em Let $V$ be a nonlocal vertex algebra. A subset  $U$ 
of $\Hom (V,V\otimes \C((x)))$ is said to be {\em $\Delta$-closed}
if for any $a(x)\in U$, there exist $a_{(1i)}(x), a_{(2i)}(x)\in U$
for $i=1,\dots,r$  such that
\begin{eqnarray}
a(x_{1})Y(v,x_{2})
=\sum_{i=1}^{r} Y(a_{(1i)}(x_{1}-x_{2})v,x_{2})a_{(2i)}(x_{1})
\;\;\;\mbox{ for all }v\in V. 
\end{eqnarray}}
\ed

In terms of this notion, for any $a(x)\in {\rm PDer}^{-}(V)$,
$\{ a(x),1\}$ is $\Delta$-closed and for 
$\psi(x)\in {\rm PEnd}^{-}(V)$, $\{ \psi(x)\}$ is $\Delta$-closed.

Note that the space $\Hom_{\C}(V,V\otimes \C((x)))$,
 which can be naturally identified with
$\End_{\C((x))}(V\otimes \C((x)))$,
is an associative algebra 
over $\C((x))$ with the formal differential operator 
$\partial=\frac{d}{dx}$ as a $\C$-linear derivation. 
 
The followings are straightforward consequences:

\bl{lfacts-delta-closedsubsets}
Let $V$ be a nonlocal vertex algebra. 
a) The linear span of any $\Delta$-closed 
subset of $\Hom (V,V\otimes \C((x)))$ is $\Delta$-closed. 
b)  The sum of $\Delta$-closed subspaces
of $\Hom (V,V\otimes \C((x)))$ is $\Delta$-closed. 
c) If $U$ is a $\Delta$-closed 
subset of $\Hom (V,V\otimes \C((x)))$, then
$U\cdot U$ (the linear span of $a(x)b(x)$ for $a(x),b(x)\in U$)
is $\Delta$-closed.
d) If $U$ is a $\Delta$-closed subspace, $U+\partial (U)$ is $\Delta$-closed.
\el

\bd{dbv}
{\em Let $V$ be a nonlocal vertex algebra.
Define $B(V)$ to be the sum of all the $\Delta$-closed subspaces $U$ of
$\Hom (V,V\otimes \C((x)))$ such that
\begin{eqnarray}
a(x){\bf 1}\in \C {\bf 1}\;\;\;\mbox{ for }a(x)\in U.
\end{eqnarray}}
\ed

\bp{pbv-algebra}
For any nonlocal vertex algebra $V$, 
$B(V)$ is a $\Delta$-closed associative subalgebra 
of $\Hom (V,V\otimes \C((x)))$ and it is closed 
under the derivation $\partial=\frac{d}{dx}$. 
Furthermore, $V$ is a module for $(B(V),\partial)$
with $Y_{V}(a(x),x_{0})=a(x_{0})$ for $a(x)\in B(V)$.
\ep

\begin{proof} From Lemma \ref{lfacts-delta-closedsubsets} b),
$B(V)$ is $\Delta$-closed, and then by
Lemma \ref{lfacts-delta-closedsubsets} c),
$B(V)\cdot B(V)$ is $\Delta$-closed. 
Clearly, we have $B(V){\bf 1}\subset \C {\bf 1}$. Then
$B(V)B(V){\bf 1}\subset \C {\bf 1}$. 
Now, it follows that
$B(V)\cdot B(V)\subset B(V)$. This proves that $B(V)$ is a subalgebra.
By Lemma \ref{lfacts-delta-closedsubsets} d),
$B(V)+\partial B(V)$ is $\Delta$-closed and
we have $\partial B(V){\bf 1}=0$.
Then $\partial B(V)\subset B(V)$.

Next, we prove that $V$ is a module for $B(V)$ viewed as
a nonlocal vertex algebra. 
First, by definition we have 
$$Y_{V}(a(x),x_{0})v=a(x_{0})v\in V\otimes \C((x_{0}))\subset V((x_{0}))
\;\;\;\mbox{ for }a(x)\in B(V),\; v\in V.$$
For $a(x),b(x)\in B(V),\; v\in V$, we have
$$Y_{V}(a(x),x_{0}+x_{2})Y_{V}(b(x),x_{2})v=a(x_{0}+x_{2})b(x_{2})v$$ and
$$Y_{V}(Y(a(x),x_{0})b(x),x_{2})v
=(a(x+x_{0})b(x)v)|_{x=x_{2}}=a(x_{2}+x_{0})b(x_{2})v,$$
noticing that $Y(a(x),x_{0})b(x)
=(e^{x_{0}\frac{d}{dx}}a(x))b(x)=a(x+x_{0})b(x)$.
Since $b(x_{2})v\in V\otimes \C(((x_{2}))$,
there exists $l\in \N$ such that
$x_{1}^{l}a(x_{1})b(x_{2})v\in V\otimes \C((x_{2}))[[x_{1}]]$. Then
$$(x_{0}+x_{2})^{l}a(x_{0}+x_{2})b(x_{2})v
=(x_{2}+x_{0})^{l}a(x_{2}+x_{0})b(x_{2})v.$$
Consequently, we have
$$(x_{0}+x_{2})^{l}Y(a(x),x_{0}+x_{2})Y(b(x),x_{2})v
=(x_{0}+x_{2})^{l}Y(Y(a(x),x_{0})b(x),x_{2})v.$$
We also have $Y_{V}(1,x)=1_{V}$.
Thus $V$ is a module for $B(V)$ viewed as a nonlocal vertex algebra.
\end{proof}

\br{recall-bialgebra}
{\em Recall that a {\em coalgebra} over $\C$ is 
a vector space $A$ equipped with 
two linear maps $\Delta: A\rightarrow A\otimes A$, 
called the {\em comultiplication map}, and
$\varepsilon: A\rightarrow \C$, called the {\em counit map}, such that
\begin{eqnarray}
& &(1\otimes \Delta) \Delta=(\Delta\otimes 1)\Delta,\label{ecoalgebra-delta}\\
& &\sum \varepsilon (a_{(1)})a_{(2)}=\sum \varepsilon (a_{(2)})a_{(1)}=a
\;\;\;\mbox{ for }a\in A.
\end{eqnarray}
A {\em bialgebra} is a (unital) associative algebra equipped with 
a coalgebra structure such that the comultiplication map
and the counit map are homomorphisms of
associative algebras.}
\er

For convenience, we formulate the following notion:

\bd{ddiff-bilagbera}
{\em A {\em differential bialgebra} is 
a bialgebra $(B,\Delta,\varepsilon)$ equipped with 
a derivation $\partial$ such that $\varepsilon \circ \partial =0$ and
$\Delta \partial =(\partial\otimes 1+1\otimes \partial)\Delta$.}
\ed

We shall prove that if 
a nonlocal vertex algebra $V$ satisfies a certain condition,
$B(V)$ has a natural bialgebra structure. 
Recall from \cite{ek} that a nonlocal vertex algebra $V$ 
is {\em nondegenerate} if for any positive integer $n$, 
the linear map
$$Z_{n}: V^{\otimes n}\otimes \C((x_{1}))\cdots ((x_{n}))
\rightarrow V((x_{1}))\cdots ((x_{n}))$$
defined by
$$Z_{n}(v_{(1)}\otimes \cdots \otimes v_{(n)}\otimes f)
=f Y(v_{(1)},x_{1})\cdots Y(v_{(n)},x_{n}){\bf 1}$$
is injective.

\br{rdef-X}
{\em Let $V$ be a nonlocal vertex algebra. For any positive integer
 $n$, we define a linear map
$$X_{n}: V^{(n+1)}\otimes \C((x_{1}))\cdots ((x_{n}))
\rightarrow V((x_{1}))\cdots ((x_{n}))$$
by
$$X_{n}(v_{(1)}\otimes \cdots \otimes v_{(n)}\otimes v\otimes f)
=f Y(v_{(1)},x_{1})\cdots Y(v_{(n)},x_{n})v.$$
It was proved in \cite{li-simple} that if $Z_{n+1}$ is injective, so is $X_{n}$. 
Thus if $V$ is nondegenerate, for every positive integer $n$,
 $X_{n}$ is injective,
in particular, $X_{1}$ and $X_{2}$ are injective.}
\er 

\bt{tbialgebra-structure}
Let $V$ be a nondegenerate nonlocal vertex algebra.
Then $B(V)$ is a differential bialgebra 
with the comultiplication 
$\Delta$ and the counit $\varepsilon$ uniquely determined by
\begin{eqnarray}
& &a(x){\bf 1}=\varepsilon (a(x)){\bf 1},\\
& &\Delta(a(x))=\sum a_{(1)}(x)\otimes a_{(2)}(x)
\end{eqnarray}
for $a(x)\in B(V)$, where 
$a(x_{1})Y(v,x_{2})=\sum
Y(a_{(1)}(x_{1}-x_{2})v,x_{2})a_{(2)}(x_{1})$ for all $v\in V$ and
we are using the coalgebra Sigma notation.
Furthermore, we have
\begin{eqnarray}
& &{\rm PDer}^{-}(V)
=\{ a\in B(V)\;|\; \Delta(a)=a\otimes 1+1\otimes a\},
\label{e3.5}\\
& &{\rm PEnd}^{-}(V)
=\{ a\in B(V)\;|\; \Delta(a)=a\otimes a\}.\label{e3.6}
\end{eqnarray}
\et

\begin{proof} First we show that $\Delta$ is well defined. It suffices to prove that if
\begin{eqnarray}\label{eif}
0\ne \sum_{r=1}^{k}a_{r}(x)\otimes b_{r}(x)\in B(V)\otimes B(V),
\end{eqnarray}
then there exist $u,v\in V$ such that
\begin{eqnarray}\label{ethen}
\sum_{r=1}^{k}Y(a_{r}(x_{1}-x_{2})u,x_{2})b_{r}(x_{1})v\ne 0.
\end{eqnarray}
{}From (\ref{eif}), there exist $u,v\in V$ such that
$$0\ne \sum_{r=1}^{k}a_{r}(-x_{2})u\otimes b_{r}(x_{1})v
\in (V\otimes \C((x_{2})))\otimes (V\otimes \C((x_{1})))
\subset (V\otimes V\otimes \C((x_{2})))((x_{1})),$$
which (by applying $e^{-x_{1}\partial/\partial x_{2}}$) implies
$$0\ne \sum_{r=1}^{k}a_{r}(-x_{2}+x_{1})u\otimes b_{r}(x_{1})v
\in (V\otimes V\otimes \C((x_{2})))((x_{1})).$$
Since the linear map $X_{1}$ is injective, we have
\begin{eqnarray}\label{enon-zero}
0\ne \sum_{r=1}^{k}Y(a_{r}(-x_{2}+x_{1})u,x_{2})b_{r}(x_{1})v
\in V((x_{2}))((x_{1})).
\end{eqnarray}
Let $l\in \N$ be such that $a_{r}(x)u\in V\otimes \C[[x]]$ for $r=1,\dots,k$. Consequently,
$$(x_{1}-x_{2})^{l}a_{r}(x_{1}-x_{2})u=(x_{1}-x_{2})^{l}a_{r}(-x_{2}+x_{1})u.$$
Multiplying (\ref{enon-zero}) by $(x_{1}-x_{2})^{l}$ and then using cancellation 
we get (\ref{ethen}).
This proves that $\Delta$ is well defined. 

Next we establish the coassociativity. Let $a(x)\in B(V)$ and let $u,v,w\in V$. 
We have
\begin{eqnarray}\label{ewehave}
& &a(x_{1})Y(u,x_{0}+x_{2})Y(v,x_{2})w\nonumber\\
&=&\sum Y(a_{(1)}(x_{1}-x_{0}-x_{2})u,x_{0}+x_{2})
a_{(2)}(x_{1})Y(v,x_{2})w\nonumber\\
&=&\sum Y(a_{(1)}(x_{1}-x_{0}-x_{2})u,x_{0}+x_{2})
Y(a_{(2,1)}(x_{1}-x_{2})v,x_{2})a_{(2,2)}(x_{1})w.\ \ \ 
\end{eqnarray}
On the other hand, we have
\begin{eqnarray}\label{eotherhand}
& &a(x_{1})Y(Y(u,x_{0})v,x_{2})w\nonumber\\
&=&\sum Y(a_{(1)}(x_{1}-x_{2})
Y(u,x_{0})v, x_{2}) a_{(2)}(x_{1})w\nonumber\\
&=&\sum 
Y(Y(a_{(1,1)}(x_{1}-x_{2}-x_{0})u,x_{0})a_{(1,2)}(x_{1}-x_{2})v,x_{2})
a_{(2)}(x_{1})w.
\end{eqnarray}
With $a_{(1,1)}(x)u,\; a_{(1,2)}(x)v,\; a_{(2)}(x)w\in V\otimes \C((x))$,
there exists $l\in \N$ such that
\begin{eqnarray*}
& &(x_{0}+x_{2})^{l}Y(Y(u,x_{0})v,x_{2})w=Y(u,x_{0}+x_{2})Y(v,x_{2})w,\\
& &(x_{0}+x_{2})^{l}\sum 
Y(Y(a_{(1,1)}(x_{1}-x_{2}-x_{0})u,x_{0})a_{(1,2)}(x_{1}-x_{2})v,x_{2})
a_{(2)}(x_{1})w\nonumber\\
& &\ \ \ \ =(x_{0}+x_{2})^{l}\sum 
Y(a_{(1,1)}(x_{1}-x_{2}-x_{0})u,x_{0}+x_{2})Y(a_{(1,2)}(x_{1}-x_{2})v,x_{2})
a_{(2)}(x_{1})w.
\end{eqnarray*}
Then using (\ref{ewehave}) and (\ref{eotherhand}) we obtain
\begin{eqnarray*}
& &\sum (x_{0}+x_{2})^{l}Y(a_{(1)}(x_{1}-x_{0}-x_{2})u,x_{0}+x_{2})
Y(a_{(2,1)}(x_{1}-x_{2})v,x_{2})a_{(2,2)}(x_{1})w\nonumber\\
&=&\sum (x_{0}+x_{2})^{l}Y(a_{(1,1)}(x_{1}-x_{2}-x_{0})u,x_{0}+x_{2})
Y(a_{(1,2)}(x_{1}-x_{2})v,x_{2})a_{(2)}(x_{1})w,\ \ \ \ 
\end{eqnarray*}
which by cancellation and substitution yields
\begin{eqnarray}
& &\sum Y(a_{(1)}(x_{1}-x_{3})u,x_{3})
Y(a_{(2,1)}(x_{1}-x_{2})v,x_{2})a_{(2,2)}(x_{1})w\nonumber\\
&=&\sum Y(a_{(1,1)}(x_{1}-x_{3})u,x_{3})
Y(a_{(1,2)}(x_{1}-x_{2})v,x_{2})a_{(2)}(x_{1})w.
\end{eqnarray}
By multiplying by certain nonnegative powers of
$(x_{1}-x_{3})$ and $(x_{2}-x_{3})$ and then using cancellation we get
\begin{eqnarray}
& &\sum Y(a_{(1)}(-x_{3}+x_{1})u,x_{3})
Y(a_{(2,1)}(-x_{2}+x_{1})v,x_{2})a_{(2,2)}(x_{1})w\nonumber\\
&=&\sum Y(a_{(1,1)}(-x_{3}+x_{1})u,x_{3})
Y(a_{(1,2)}(-x_{2}+x_{1})v,x_{2})a_{(2)}(x_{1})w.
\end{eqnarray}
Since $X_{2}$ is injective, we have
$$\sum a_{(1)}(x_{3})u\otimes a_{(2,1)}(x_{2})v\otimes a_{(2,2)}(x_{1})w
=\sum a_{(1,1)}(x_{3})u\otimes a_{(1,2)}(x_{2})\otimes a_{(2)}(x_{1})w. $$
This gives the coassociativity. It is straightforward to show that
%For $a(x),b(x)\in B(V),\; v\in V$, we have
%\begin{eqnarray*}
%a(x_{1})b(x_{1})Y(v,x_{2})
%&=&\sum a(x_{1})Y(b_{(1)}(x_{1}-x_{2})v,x_{2})b_{(2)}(x_{1}) \\
%&=&\sum Y(a_{(1)}(x_{1}-x_{2})b_{(1)}(x_{1}-x_{2})v,x_{2})a_{(2)}(x_{1})b_{(2)}(x_{1}),
%\end{eqnarray*}
%which gives $\Delta (a(x)b(x))=\Delta(a(x))\Delta(b(x))$.
%We also have $\Delta({\bf 1})={\bf 1}\otimes {\bf 1}$.
$\Delta$ is an algebra homomorphism.

As  $a(x){\bf 1}\in \C {\bf 1}$ for $a(x)\in B(V)$ by assumption, 
the map $\varepsilon: B(V) \rightarrow \C$ is well defined. 
It is straightforward to see that $\varepsilon$ is an algebra homomorphism.
For $a(x)\in B(V)$, we have
\begin{eqnarray}
a(x_{1})&=&a(x_{1})Y({\bf 1},x_{2})\nonumber\\
&=&\sum Y(a_{(1)}(x_{1}-x_{2}){\bf 1},x_{2})a_{(2)}(x_{1})\nonumber\\
&=&\sum Y(\varepsilon(a_{(1)}(x)){\bf 1},x_{2})a_{(2)}(x_{1})\nonumber\\
&=&\sum \varepsilon(a_{(1)}(x))a_{(2)}(x_{1})
\end{eqnarray}
and
\begin{eqnarray}
a(x)v&=&\Res_{x_{2}}x_{2}^{-1}a(x_{1})Y(v,x_{2}){\bf 1}\nonumber\\
&=&\sum \Res_{x_{2}}x_{2}^{-1}
Y(a_{(1)}(x_{1}-x_{2})v,x_{2})a_{(2)}(x_{1}){\bf 1}\nonumber\\
&=&\sum \Res_{x_{2}}x_{2}^{-1}\varepsilon(a_{(2)}(x))
Y(a_{(1)}(x_{1}-x_{2})v,x_{2}){\bf 1}\nonumber\\
&=&\sum \Res_{x_{2}}x_{2}^{-1}
\varepsilon(a_{(2)}(x))e^{-x_{2}\frac{\partial}{\partial x_{1}}}
Y(a_{(1)}(x_{1})v,x_{2}){\bf 1}\nonumber\\
&=&\sum \varepsilon(a_{(2)}(x))a_{(1)}(x_{1})v.
\end{eqnarray}
This proves that $\varepsilon$ is a counit.
Therefore, $(B(V), \Delta, \varepsilon)$ is a bialgebra.

It is straightforward to show 
$$\Delta \partial =(\partial\otimes 1+1\otimes \partial)\Delta,\;\;
\;\; \varepsilon\circ \partial=0. $$
%Let $a(x)\in B(V)$. For $v\in V$, we have
%\begin{eqnarray*}
%a'(x)Y(v,x_{2})
%&=&\sum \frac{\partial}{\partial x}Y(a_{(1)}(x-x_{2})v,x_{2})a_{(2)}(x)\\
%&=&\sum Y(a'_{(1)}(x-x_{2})v,x_{2})a_{(2)}(x)+
%Y(a_{(1)}(x-x_{2})v,x_{2})a'_{(2)}(x).
%\end{eqnarray*}
%It follows that
%$$\Delta \partial a(x)=(\partial\otimes 1+1\otimes \partial)\Delta (a(x)).$$
%For $a(x)\in B(V)$, we also have
%$\varepsilon(\partial a(x)){\bf 1}=(\partial a(x)){\bf 1}
%=\frac{d}{dx}a(x){\bf 1}=0$.
Thus $B(V)$ is a differential bialgebra.
The assertions (\ref{e3.5}) and (\ref{e3.6}) are clear.
\end{proof}

We have the following consequence:

\bc{cmain-bialgebra}
Let $V$ be a nonlocal vertex algebra. 
Assume that $V$ is of countable dimension over $\C$ and that
$V$ as a $V$-module is irreducible.
Then all the assertions of Theorem \ref{tbialgebra-structure} hold.
In particular, if $V$ is a simple vertex operator algebra in the sense 
of \cite{flm}, all the assertions of Theorem \ref{tbialgebra-structure} hold.
\ec

\begin{proof} Under the assumptions, by \cite{li-simple} and
\cite{li-qva2}, $V$ is nondegenerate.
Then it follows immediately from Theorem \ref{tbialgebra-structure}.
\end{proof}

%\bc{cmain-bialgebra-pbw}
%Let $V$ be a nonlocal vertex algebra. Suppose that $V$ 
%is of countable dimension over $\C$  and that $V$ has a P-B-W basis.
%Then all the assertions of Theorem \ref{tbialgebra-structure} hold.
%\ec

%\begin{proof}
%It was proved in \cite{li-qva2} that $V$ is nondegenerate.
%\end{proof}

\bp{pcocommutativity}
Let $V$ be a nondegenerate ordinary vertex algebra and let
$B(V)^{c}$ be the sum of all the $\Delta$-closed subspaces $U$ of $\Hom (V,V\otimes \C((x)))$ 
such that
\begin{eqnarray}
& &a(x){\bf 1}\in \C {\bf 1},\\
& &[\D,a(x)]=\frac{d}{dx}a(x)\;\;\;\mbox{ for }a(x)\in U.\label{e3.15}
\end{eqnarray}
Then $B(V)^{c}$ is a cocommutative sub-bialgebra of $B(V)$.
\ep

\begin{proof} It is straightforward to show that $B(V)^{c}$ 
is a sub-differential bialgebra of $B(V)$. Now we prove the cocommutativity.
Let $a(x)\in B(V)^{c}$ with $\Delta(a(x))=\sum a_{(1)}(x)\otimes a_{(2)}(x)$. 
For $u,v\in V$, we have
\begin{eqnarray*}
a(x_{1})Y(u,x)v=\sum Y(a_{(1)}(x_{1}-x)u,x)a_{(2)}(x_{1})v.
\end{eqnarray*}
Using the skew symmetry of the vertex algebra $V$ and (\ref{e3.15}) we get
\begin{eqnarray*}
& &a(x_{1})Y(u,x)v=a(x_{1})e^{x\D}Y(v,-x)u=e^{x\D}a(x_{1}-x)Y(v,-x)u,\\
& &\sum Y(a_{(1)}(x_{1}-x)u,x)a_{(2)}(x_{1})v=\sum e^{x\D}Y(a_{(2)}(x_{1})v,-x)a_{(1)}(x_{1}-x)u.
\end{eqnarray*}
Consequently, 
$$a(x_{1}-x)Y(v,-x)u=\sum Y(a_{(2)}(x_{1})v,-x)a_{(1)}(x_{1}-x)u,$$
which is,
$$a(x_{1})Y(v,x)u=\sum Y(a_{(2)}(x_{1}-x)v,x)a_{(1)}(x_{1})u.$$
Thus $\Delta(a(x))=\sum a_{(2)}(x)\otimes a_{(1)}(x)$. This proves the cocommutativity.
\end{proof}

\section{Smash product nonlocal vertex algebras} 

In this section we formulate a notion of vertex bialgebra and a notion
of module nonlocal vertex algebra for a vertex bialgebra.  We then
give a smash product construction of nonlocal vertex algebras.

The following is our notion of vertex bialgebra:

\bd{dvertex-bialgebra}
{\em A {\em vertex bialgebra} is a nonlocal vertex algebra $V$ 
equipped with a coalgebra structure $(V,\Delta,\varepsilon)$ such that 
$\Delta: V\rightarrow V\otimes V$ and $\varepsilon: V\rightarrow \C$
are  homomorphisms of nonlocal vertex algebras.}
\ed

\bex{exadiff-bialgebra-examples}
{\em Here we extend Borcherds' construction of vertex algebras for
vertex bialgebras. 
Let $(B,\Delta,\varepsilon,\partial)$ be a differential bialgebra
(see Definition \ref{ddiff-bilagbera}).
First we have a nonlocal vertex algebra $(B,\partial)$.
{}From definition we have $\varepsilon(1)=1$ and $\Delta(1)=1\otimes 1$.
Furthermore, for $a,b\in B$ we have
\begin{eqnarray*}
& &\varepsilon (Y(a,x)b)=\varepsilon ((e^{x\partial}a)b)
=\varepsilon (e^{x\partial}a)\varepsilon (b)=\varepsilon(a)\varepsilon(b),\\
& &\Delta (Y(a,x)b)=\Delta ((e^{x\partial}a)b)
=\Delta (e^{x\partial}a)\Delta(b)
=(e^{x(\partial\otimes 1+1\otimes\partial)}\Delta(a))\Delta(b)
=Y(\Delta (a),x)\Delta(b).
\end{eqnarray*}
This shows that $\Delta$ and $\varepsilon$ are homomorphisms of 
nonlocal vertex algebras. Therefore, $B$ is a vertex bialgebra.}
\eex

\bex{exahopf-Lie}
{\em Let $\g$ be a Lie algebra equipped with a derivation $d$.
Then $U(\g)$ is a Hopf algebra with $d$ naturally extended to a derivation
of $U(\g)$. It is straightforward to show that
$U(\g)$ viewed as a bialgebra equipped with the derivation $d$ is a differential
bialgebra. In view of Example \ref{exadiff-bialgebra-examples}, 
$U(\g)$ is naturally a vertex bialgebra.}
\eex

\bex{exa-vla}
{\em Here we show that vertex algebras associated with vertex Lie algebras
can be made into vertex bialgebras.
Let $C$ be a vertex Lie algebra (see [P]), or namely a conformal Lie algebra
(see [K]). Associated to $C$ we have an honest Lie algebra
$\L_{C}$. As a vector space, 
$$\L_{C}=(C\otimes \C[t,t^{-1}])/(\partial\otimes 1+1\otimes d/dt)(C\otimes \C[t,t^{-1}]).$$
For $a\in C,\; n\in \Z$, we denote by $a(n)$ the element of $\L_{C}$ 
corresponding to $a\otimes t^{n}$.
The Lie algebra $\L_{C}$ has a polar decomposition $\L_{C}=\L_{C}^{+}\oplus \L_{C}^{-}$ 
into subalgebras, where
$\L_{C}^{+}$ is spanned by $a(n)$ for $a\in C,\; n\ge 0$ and
$\L_{C}^{-}$ is spanned by $a(n)$ for $a\in C,\; n< 0$.
Form the induced module of $\L_{C}$ from the trivial $\L_{C}^{+}$-module $\C$:
$$V_{C}=U(\L_{C})\otimes _{U(\L_{C}^{+})}\C.$$ 
In view of the P-B-W theorem, $V_{C}=U(\L_{C}^{-})$ as vector spaces.
Set 
$${\bf 1}=1\otimes 1\in V_{C}.$$
Identify $C$ as a subspace of $V_{C}$ through the linear map
$u\mapsto u(-1){\bf 1}$.
Then (see [DLM], [P]) there exists a unique vertex algebra structure on 
$V_{C}$ with ${\bf 1}$ as the vacuum vector and with
$$Y(u,x)=u(x)=\sum_{n\in \Z}u(n)x^{-n-1}\;\;\;\mbox{ for }u\in C.$$
Equip $V_{C}\otimes V_{C}$ with the tensor product $\L_{C}$-module structure.
We have $\L_{C}^{+}({\bf 1}\otimes {\bf 1})=0$.
Then there exists a unique $\L_{C}$-module map $\Delta$ from $V_{C}$ to $V_{C}\otimes V_{C}$ 
such that $\Delta({\bf 1})={\bf 1}\otimes {\bf 1}$. For $a\in C$, we have
$$\Delta(a)=\Delta (a(-1){\bf 1})=a(-1)\Delta ({\bf 1})
=a(-1){\bf 1}\otimes {\bf 1}+{\bf 1}\otimes a(-1){\bf 1}=a\otimes {\bf 1}+{\bf 1}\otimes a.$$
Then
$$Y(\Delta(a),x)=Y(a,x)\otimes 1+1\otimes Y(a,x)=a(x)\otimes 1+1\otimes a(x).$$
Now, for any $a\in C,\; v\in V_{C}$, we have 
$$\Delta(Y(a,x)v)=\Delta(a(x)v)=a(x)\Delta(v)=
(a(x)\otimes 1+1\otimes a(x))\Delta(v)=Y(\Delta(a),x)\Delta(v).$$
As $C$ generates $V_{C}$ as a vertex algebra, $\Delta$ is a vertex algebra homomorphism.
Similarly, equipping $\C$ with the trivial $\L_{C}$-module, we have
a unique $\L_{C}$-module map $\varepsilon: V_{C}\rightarrow \C$ 
such that  $\varepsilon({\bf 1})=1$. For $a\in C,\; v\in V_{C}$, as
$\varepsilon(a)=\varepsilon(a(-1){\bf 1})=a(-1)\varepsilon({\bf 1})=0$, we have
$$\varepsilon (Y(a,x)v)=\varepsilon (a(x)v)=a(x)\varepsilon (v)=0
=Y(\varepsilon(a),x)\varepsilon(v).$$
{}From the same reasoning $\varepsilon$ is a vertex algebra homomorphism.
Therefore, the associated vertex algebra $V_{C}$ equipped with
the maps $\Delta$ and $\varepsilon$ defined above is a vertex bialgebra.}
\eex

The following is a vertex analogue of
the notion of module algebra for a bialgebra (cf. \cite{mon}):

\bd{dmodule-algebra}
{\em Let $H$ be a vertex bialgebra.
An {\em $H$-module nonlocal vertex algebra} is a nonlocal vertex algebra $V$ 
equipped with an $H$-module structure such that 
\begin{eqnarray}
& &Y(h,x)v\in V\otimes \C((x)),\\
& &Y(h,x){\bf 1}_{V}=\varepsilon(h){\bf 1}_{V},\\
& &Y(h,x_{1})Y(u,x_{2})v
=\sum Y(Y(h_{(1)},x_{1}-x_{2})u,x_{2})Y(h_{(2)},x_{1})v
\end{eqnarray}
for $h\in H,\; u,v\in V$, where ${\bf 1}_{V}$ denotes the vacuum vector of $V$.
}
\ed

The following result is very useful in constructing module nonlocal vertex algebras:

\bl{lmodule-algebra}
Let $H$ be a vertex bialgebra and let $V$ be a nonlocal vertex algebra.
Suppose that $(V,Y_{V}^{H})$ is an $H$-module and
$S$ is a generating subset of $H$ as a nonlocal vertex algebra
such that for $h\in S,\;u, v\in V$,
\begin{eqnarray}
& &Y_{V}^{H}(h,x)\in \Hom (V,V\otimes \C((x))),\\
& &Y_{V}^{H}(h,x){\bf 1}=\varepsilon(h){\bf 1},\label{ecounit-generatingr}\\
& &Y_{V}^{H}(h,x_{1})Y(u,x_{2})v
=\sum Y(Y_{V}^{H}(h_{(1)},x_{1}-x_{2})u,x_{2})Y_{V}^{H}(h_{(2)}),x_{1})v.
\label{ema-generating}
\end{eqnarray}
Then $V$ equipped with the $H$-module structure $Y_{V}^{H}$
is an $H$-module nonlocal vertex algebra.
\el

\begin{proof} Let $g,h\in H$ be such that
$$Y_{V}^{H}(g,x),\;\; Y_{V}^{H}(h,x)\in \Hom (V,V\otimes \C((x))).$$
For any $v\in V$, as $(V,Y_{V}^{H})$ is an $H$-module, 
there exists a nonnegative integer $l$ such that
\begin{eqnarray}\label{eweak-assoc-proof}
(x_{0}+x_{2})^{l}Y_{V}^{H}(Y(g,x_{0})h,x_{2})v
=(x_{0}+x_{2})^{l}Y_{V}^{H}(g,x_{0}+x_{2})Y_{V}^{H}(h,x_{2})v.
\end{eqnarray}
Since $Y_{V}^{H}(h,x_{2})v\in V\otimes \C((x_{2}))$, there exists $l'\in \N$
such that $x^{l'}Y_{V}^{H}(g,x)Y_{V}^{H}(h,x_{2})v\in V[[x]]\otimes \C((x_{2}))$.
Then
$$(x_{0}+x_{2})^{l'}Y_{V}^{H}(g,x_{0}+x_{2})Y_{V}^{H}(h,x_{2})v
=(x_{2}+x_{0})^{l'}Y_{V}^{H}(g,x_{2}+x_{0})Y_{V}^{H}(h,x_{2})v.$$
Combining this with (\ref{eweak-assoc-proof}) we get
$$(x_{0}+x_{2})^{l+l'}Y_{V}^{H}(Y(g,x_{0})h,x_{2})v
=(x_{2}+x_{0})^{l+l'}Y_{V}^{H}(g,x_{2}+x_{0})Y_{V}^{H}(h,x_{2})v.$$
By cancellation we get
$$Y_{V}^{H}(Y(g,x_{0})h,x_{2})v
=Y_{V}^{H}(g,x_{2}+x_{0})Y_{V}^{H}(h,x_{2})v
=\left((e^{x_{0}\frac{d}{dx}}Y_{V}^{H}(g,x))Y_{V}^{H}(h,x)\right)|_{x=x_{2}}v.$$
Since $S$ generates $H$, we have $Y_{V}^{H}(h,x)\in \Hom (V,V\otimes \C((x)))$ for all $h\in H$.
Because $\Delta,\varepsilon$ are homomorphisms of nonlocal vertex algebras and
because $S$ generates $H$ as a nonlocal vertex  algebra, it follows that
(\ref{ecounit-generatingr}) and (\ref{ema-generating}) hold for all $h\in H$.
Thus $V$ is an $H$-module nonlocal vertex algebra.
\end{proof}

We have the following simple result:

\bl{lmodule-algebra-consequence}
Let $H$ be a vertex bialgebra and let $V$ be an $H$-module nonlocal vertex algebra
with the $H$-module structure $Y_{V}^{H}$.
Then $Y_{V}^{H}$ is a homomorphism of nonlocal vertex algebras from $H$ to $B(V)$.
In particular,
\begin{eqnarray}
Y_{V}^{H}(Y(g,x_{0})h,x)v=\left(e^{x_{0}\frac{d}{dx}}Y_{V}^{H}(g,x)\right)Y_{V}^{H}(h,x)v
\end{eqnarray}
for $g,h\in H,\; v\in V$. 
\el

\begin{proof} {}From Definition \ref{dmodule-algebra}, we have
$Y_{V}^{H}(h,x)\in \Hom(V,V\otimes \C((x)))$,
 $Y_{V}^{H}(h,x){\bf 1}=\varepsilon(h){\bf 1}\in \C {\bf 1}$ for $h\in H$
and $\{ Y_{V}^{H}(h,x)\;|\; h\in H\}$ is a $\Delta$-closed subspace.
Then the linear map $Y_{V}^{H}$ maps $H$ into $B(V)$. From the proof of
Lemma \ref{lmodule-algebra}, the particular assertion holds. It follows that
$Y_{V}^{H}$ is a homomorphism of nonlocal vertex algebras from $H$ to $B(V)$.
\end{proof}

Recall from Corollary \ref{cmain-bialgebra} that if
$V$ is nondegenerate, then $B(V)$ is a differential bialgebra, 
hence a vertex bialgebra.

\bp{pbialgebra}
Let $V$ be a nondegenerate nonlocal vertex algebra.
Then $V$ is a $B(V)$-module nonlocal vertex algebra
with $Y(a(x),x_{0})=a(x_{0})$ for $a(x)\in B(V)$.
Furthermore, for any vertex bialgebra $H$, 
an $H$-module nonlocal vertex algebra structure $Y_{V}^{H}$
amounts to a homomorphism  of vertex bialgebras from $H$ to $B(V)$.
\ep

\begin{proof} By Proposition \ref{pbv-algebra},
$V$ is a module for $B(V)$ viewed as a nonlocal vertex algebra.
For $a(x)\in B(V),\; v\in V$, we have
$$Y(a(x),x_{0})v=a(x_{0})v\in V\otimes \C((x_{0})),$$
and
\begin{eqnarray*}
& &Y(a(x),x_{0}){\bf 1}=a(x_{0}){\bf 1}=\varepsilon(a(x)){\bf 1},\\
& &Y(a(x),x_{1})Y(v,x_{2})=a(x_{1})Y(v,x_{2})
=\sum Y(a_{(1)}(x_{1}-x_{2})v,x_{2})a_{(2)}(x_{1})
\\
& &\ \ \ \ \ \ \ \ =\sum Y(Y(a_{(1)}(x),x_{1}-x_{2})v,x_{2})Y(a_{(2)}(x),x_{1}).
\end{eqnarray*}
Therefore $V$ is a $B(V)$-module nonlocal vertex algebra.

Let $H$ be a vertex bialgebra and  assume that $(V,Y_{V}^{H})$ is an $H$-module nonlocal vertex
algebra.  {}From Lemma \ref{lmodule-algebra-consequence},
$Y_{V}^{H}$ is a homomorphism of nonlocal vertex algebras from $H$ to $B(V)$. 
For $h\in H$, since $Y_{V}^{H}(h,x){\bf 1}=\varepsilon(h){\bf 1}$, by definition we have
$\varepsilon (Y(h,x))=\varepsilon(h)$. For $h\in H,\; v\in V$, we have
$$Y(h,x_{1})Y(v,x)=\sum Y(Y(h_{(1)},x_{1}-x)v,x)Y(h_{(2)},x_{1}).$$
By definition we have
$$\Delta (Y(h,x))=\sum Y(h_{(1)},x)\otimes Y(h_{(2)},x)\in B(V)\otimes B(V).$$
Therefore $Y_{V}^{H}$ is a homomorphism of vertex bialgebras.

On the other hand, assume that $Y_{V}^{H}$ is a homomorphism
of vertex bialgebras from $H$ to $B(V)$. As $V$ is a $B(V)$-module
nonlocal vertex algebra, it is straightforward to see that 
$V$ equipped with $Y_{V}^{H}$ is an $H$-module nonlocal vertex algebra.
\end{proof}

The following is our smash product construction of nonlocal vertex algebras:

\bt{tsmahsed-product}
Let $H$ be a vertex bialgebra and let $V$ be 
an $H$-module nonlocal vertex algebra.
Set $V\;\sharp \;H=V\otimes H$ as a vector space. 
For $u,v\in V,\; h,k\in H$, define
\begin{eqnarray}
Y^{\sharp}(u\otimes h,x)(v\otimes k)
=\sum Y(u,x)Y(h_{(1)},x)v\otimes Y(h_{(2)},x)k.
\end{eqnarray}
Then $V\;\sharp \;H$ is a nonlocal vertex algebra with $V$ and $H$ as
subalgebras and the following relation holds for $h\in H,\; u\in V$:
\begin{eqnarray}
Y^{\sharp}(h,x_{1})Y^{\sharp}(u,x_{2})=\sum
Y^{\sharp}(Y(h_{(1)},x_{1}-x_{2})u,x_{2})Y^{\sharp}(h_{(2)},x_{1}).
\end{eqnarray}
\et

\begin{proof} {}From definition, for $v\in V,\; h,k\in H$ we have
$Y(h,x)v\in V\otimes \C((x))$ and $Y(h,x)k\in H((x))$. It follows that
$$Y^{\sharp}(u\otimes h,x)(v\otimes k)\in (V\otimes H)((x))
\;\;\;\mbox{ for }u,v\in V,\; h,k\in H.$$
As $\Delta({\bf 1})={\bf 1}\otimes {\bf 1}$, 
we have $Y^{\sharp}({\bf 1}\otimes {\bf 1},x)=1$.
Using the counit property of $\varepsilon$ we also have
\begin{eqnarray*}
Y^{\sharp}(u \otimes h,x)({\bf 1}\otimes {\bf 1})
&=&\sum Y(u,x)Y(h_{(1)},x){\bf 1}\otimes Y(h_{(2)},x){\bf 1}\\
&=&\sum \varepsilon(h_{(1)})Y(u,x){\bf 1}\otimes Y(h_{(2)},x){\bf 1}\\
&=&\sum Y(u,x){\bf 1}\otimes Y(\varepsilon (h_{(1)}) h_{(2)},x){\bf 1}\\
&=&Y(u,x){\bf 1}\otimes Y(h,x){\bf 1},
\end{eqnarray*}
{}from which we have
$$Y^{\sharp}(u \otimes h,x)({\bf 1}\otimes {\bf 1})\in (V\otimes H)[[x]]
\;\;\;\mbox{ and }
\ \ \ \lim_{x\rightarrow 0}Y^{\sharp}(u \otimes h,x)({\bf 1}\otimes {\bf 1})
=u\otimes h.$$
Next we prove the weak associativity.
Let $u,v,w\in V,\; g,h,k\in H$. We have
\begin{eqnarray*}
& &Y^{\sharp}(u\otimes g,x_{0}+x_{2})
Y^{\sharp}(v\otimes h,x_{2})(w\otimes k)\\
&=&\sum Y^{\sharp}(u\otimes g,x_{0}+x_{2})
( Y(v,x_{2})Y(h_{(1)},x_{2})w\otimes Y(h_{(2)},x_{2})k)\\
&=&\sum Y(u,x_{0}+x_{2})Y(g_{(1)},x_{0}+x_{2})Y(v,x_{2})Y(h_{(1)},x_{2})w
\otimes Y(g_{(2)},x_{0}+x_{2})Y(h_{(2)},x_{2})k\\
&=&\sum Y(u,x_{0}+x_{2})Y(Y(g_{(1,1)},x_{0})v,x_{2})
Y(g_{(1,2)},x_{0}+x_{2})Y(h_{(1)},x_{2})w\\
& &\otimes Y(g_{(2)},x_{0}+x_{2})Y(h_{(2)},x_{2})k.
\end{eqnarray*}
Notice that $Y(g_{(1,2)},x_{0}+x_{2})Y(h_{(1)},x_{2})w
\in V\otimes \C((x_{0}+x_{2}))\otimes \C((x_{2}))$ and
$Y(g_{(1,1)},x_{0})v\in V\otimes \C((x_{0}))$. Then there exists $l\in \N$ such that
\begin{eqnarray*}
& &(x_{0}+x_{2})^{l}\sum Y(u,x_{0}+x_{2})Y(Y(g_{(1,1)},x_{0})v,x_{2})
Y(g_{(1,2)},x_{0}+x_{2})Y(h_{(1)},x_{2})w\\
&=&(x_{0}+x_{2})^{l}\sum Y(Y(u,x_{0})Y(g_{(1,1)},x_{0})v,x_{2})
Y(g_{(1,2)},x_{0}+x_{2})Y(h_{(1)},x_{2})w.
\end{eqnarray*}
Certainly, if necessary we can replace $l$ with a large one so that we also have
\begin{eqnarray*}
& &(x_{0}+x_{2})^{l}Y(g_{(1,2)},x_{0}+x_{2})Y(h_{(1)},x_{2})w
=(x_{0}+x_{2})^{l}Y(Y(g_{(1,2)},x_{0})h_{(1)},x_{2})w,\\
& &(x_{0}+x_{2})^{l}Y(g_{(2)},x_{0}+x_{2})Y(h_{(2)},x_{2})k
=(x_{0}+x_{2})^{l}Y(Y(g_{(2)},x_{0})h_{(2)},x_{2})k.
\end{eqnarray*}
Then we obtain
\begin{eqnarray*}
& &(x_{0}+x_{2})^{l}Y^{\sharp}(u\otimes g,x_{0}+x_{2})
Y^{\sharp}(v\otimes h,x_{2})(w\otimes k)\nonumber\\
&=&(x_{0}+x_{2})^{l}
\sum Y(Y(u,x_{0})Y(g_{(1,1)},x_{0})v,x_{2})Y(Y(g_{(1,2)},x_{0})h_{(1)},x_{2})w\nonumber\\
& &\otimes Y(Y(g_{(2)},x_{0})h_{(2)},x_{2})k.
\end{eqnarray*}
On the other hand, we have
\begin{eqnarray*}
& &Y^{\sharp}(Y^{\sharp}(u\otimes g,x_{0})(v\otimes h),x_{2})(w\otimes k)\\
&=&\sum Y^{\sharp}\left(Y(u,x_{0})Y(g_{(1)},x_{0})v
\otimes Y(g_{(2)},x_{0})h, x_{2}\right)(w\otimes k)\\
&=&\sum Y(Y(u,x_{0})Y(g_{(1)},x_{0})v,x_{2})
Y(Y(g_{(2,1)},x_{0})h_{(1)},x_{2})w
\otimes Y(Y(g_{(2,2)},x_{0})h_{(2)},x_{2})k,
\end{eqnarray*}
where we are using the property of $\Delta$ being a homomorphism of nonlocal vertex algebras,
which gives
$$\Delta(Y(g_{(2)},x)h)=Y(\Delta(g_{(2)}),x)\Delta(h)
=\sum Y(g_{(2,1)},x)h_{(1)}\otimes Y(g_{(2,2)},x)h_{(2)}.$$
By the coassociativity of $\Delta$ we have
$$\sum g_{(1,1)}\otimes g_{(1,2)}\otimes h_{(1)}\otimes g_{(2)}\otimes
h_{(2)} =\sum g_{(1)}\otimes g_{(2,1)}\otimes h_{(1)}\otimes g_{(2,2)}\otimes
h_{(2)}.$$
Putting everything together we get
\begin{eqnarray}
& &(x_{0}+x_{2})^{l}Y^{\sharp}(u\otimes g,x_{0}+x_{2})
Y^{\sharp}(v\otimes h,x_{2})(w\otimes k)\nonumber\\
&=&(x_{0}+x_{2})^{l}Y^{\sharp}(Y^{\sharp}(u\otimes g,x_{0})(v\otimes h),x_{2})(w\otimes k).
\end{eqnarray}
This proves the weak associativity.
Thus $V\;\sharp \;H$ is a nonlocal vertex algebra.

Clearly, $V$, identified with $V\otimes \C$, 
is a subalgebra of $V\;\sharp \;H$.
For $g,h\in H$, we have
$$Y^{\sharp}(g,x)({\bf 1}\otimes h)
=\sum Y(g_{(1)},x){\bf 1}\otimes Y(g_{(2)},x)h
=\sum \varepsilon (g_{(1)}) ({\bf 1}\otimes Y(g_{(2)},x)h)
={\bf 1}\otimes Y(g,x)h$$
as $g=\sum \varepsilon (g_{(1)}) g_{(2)}$. Thus $H$, identified with
$\C \otimes H$, is a subalgebra.

For $u,v\in V,\; h,k\in H$, we have
\begin{eqnarray*}
& &Y^{\sharp}(h,x_{1})Y^{\sharp}(u,x_{2})(v\otimes k)\\
&=&Y^{\sharp}(h,x_{1})(Y(u,x_{2})v\otimes k)\\
&=&\sum Y(h_{(1)},x_{1})Y(u,x_{2})v\otimes Y(h_{(2)},x_{1})k\\
&=&\sum Y(Y(h_{(1,1)},x_{1}-x_{2})u,x_{2})
Y(h_{(1,2)},x_{1})v\otimes Y(h_{(2)},x_{1})k
\end{eqnarray*}
and
\begin{eqnarray*}
& &\sum Y^{\sharp}(Y(h_{(1)},x_{1}-x_{2})u,x_{2})
Y^{\sharp}(h_{(2)},x_{1})(v\otimes k)\\
&=&\sum Y(Y(h_{(1,1)},x_{1}-x_{2})u,x_{2})Y(h_{(2,1)},x_{1})v
\otimes Y(h_{(2,2)},x_{1})k.
\end{eqnarray*}
As $(\Delta \otimes 1)\Delta (h)=(1\otimes \Delta)\Delta(h)$, we have
$$\sum h_{(1,1)}\otimes h_{(1,2)}\otimes h_{(2)}=\sum h_{(1)}\otimes
h_{(2,1)}\otimes h_{(2,2)}.$$
Therefore
\begin{eqnarray*}
Y^{\sharp}(h,x_{1})Y^{\sharp}(u,x_{2})(v\otimes k)=
\sum Y^{\sharp}(Y(h_{(1)},x_{1}-x_{2})u,x_{2})
Y^{\sharp}(h_{(2)},x_{1})(v\otimes k).
\end{eqnarray*}
Now the proof is complete.
\end{proof}

\bex{ecase-lie}
{\em Let $\g$ be a Lie algebra equipped with a derivation $d$.
Recall from Example \ref{exahopf-Lie} that the universal enveloping algebra 
$U(\g)$ is naturally a differential bialgebra and hence a vertex bialgebra. 
Let $V$ be a nonlocal vertex algebra and let $\psi: \g \rightarrow {\rm PDer}^{-}(V)$ 
be a homomorphism of differential Lie algebras.
The map $\psi$ extends uniquely to a homomorphism $\bar{\psi}: U(\g)\rightarrow B(V)$
of differential algebras. This makes $V$ a $(U(\g),d)$-module.
By Lemma \ref{lmodule-algebra}, $V$ is a $U(\g)$-module nonlocal vertex algebra.
Therefore we have the smash product nonlocal vertex algebra $V\sharp U(\g)$.}
\eex

The following result gives a construction of $V\sharp H$-modules:

\bp{pmodule-smash}
Let $H$ be a vertex bialgebra and let $V$ be an $H$-module nonlocal vertex algebra.
Suppose that $W$ is a $V$-module and an $H$-module such that
\begin{eqnarray}
& &Y(h,x)w\in W\otimes \C((x)),\label{einduction1}\\
& &Y(h,x_{1})Y(v,x_{2})w
=\sum Y(Y(h_{(1)},x_{1}-x_{2})v,x_{2})Y(h_{(2)},x_{1})w\label{einduction2}
\end{eqnarray}
for $h\in S,\; v\in V,\; w\in W$, where $S$ is a generating subset of $H$ 
as a nonlocal vertex algebra. Then $W$ is a module for $V\sharp H$ with
\begin{eqnarray}
Y_{W}(v\otimes h,x)w=Y(v,x)Y(h,x)w
\end{eqnarray}
for $h\in H,\; v\in V,\; w\in W$.
\ep

\begin{proof} The proof of Lemma \ref{lmodule-algebra} shows that
(\ref{einduction1}) and (\ref{einduction2}) hold for all $h\in H$.
For $v\in V,\; h\in H,\; w\in W$, we have
$$Y_{W}(v\otimes h,x)w=Y(v,x)Y(h,x)w\in W((x))$$
as $Y(h,x)w\in W\otimes \C((x))$ and $Y(v,x)\in \Hom(W,W((x)))$. 
{}From definition we have $Y_{W}({\bf 1}\otimes 1,x)=Y({\bf 1},x)Y(1,x)=1$.
Let $u,v\in V,\; g,h\in H,\; w\in W$. We have
\begin{eqnarray*}
& &Y_{W}(u\otimes g,x_{0}+x_{2})Y_{W}(v\otimes h,x_{2})w\\
&=&Y(u,x_{0}+x_{2})Y(g,x_{0}+x_{2})Y(v,x_{2})Y(h,x_{2})w\\
&=&\sum Y(u,x_{0}+x_{2})Y(Y(g_{(1)},x_{0})v,x_{2})
Y(g_{(2)},x_{0}+x_{2})Y(h,x_{2})w.
\end{eqnarray*}
Notice that
$$Y(g_{(2)},x_{0}+x_{2})Y(h,x_{2})w\in W\otimes \C((x_{2},(x_{0}+x_{2}))),
\ \ \ Y(g_{(1)},x_{0})v\in V\otimes \C((x_{0})).$$
Then there exists a nonnegative integer $l$ such that
\begin{eqnarray*}
& &\sum (x_{0}+x_{2})^{l}Y(u,x_{0}+x_{2})
Y(Y(g_{(1)},x_{0})v,x_{2})Y(g_{(2)},x_{0}+x_{2})Y(h,x_{2})w\\
&=&\sum (x_{0}+x_{2})^{l}Y(Y(u,x_{0})
Y(g_{(1)},x_{0})v,x_{2})Y(Y(g_{(2)},x_{0})h,x_{2})w.
\end{eqnarray*}
Thus
\begin{eqnarray*}
& &(x_{0}+x_{2})^{l}Y_{W}(u\otimes g,x_{0}+x_{2})Y_{W}(v\otimes h,x_{2})w\\
&=&\sum (x_{0}+x_{2})^{l}Y(Y(u,x_{0})Y(g_{(1)},x_{0})v,x_{2})Y(Y(g_{(2)},x_{0})h,x_{2})w.
\end{eqnarray*}
On the other hand, we have
\begin{eqnarray*}
& &Y_{W}(Y_{\sharp}(u\otimes g,x_{0})(v\otimes h),x_{2})w\\
&=&\sum Y_{W}(Y(u,x_{0})Y(g_{(1)},x_{0})v\otimes Y(g_{(2)},x_{0})h,x_{2})w\\
&=&\sum Y(Y(u,x_{0})Y(g_{(1)},x_{0})v,x_{2})Y(Y(g_{(2)},x_{0})h,x_{2})w.
\end{eqnarray*}
Therefore 
\begin{eqnarray*}
& &(x_{0}+x_{2})^{l}Y_{W}(u\otimes g,x_{0}+x_{2})Y_{W}(v\otimes h,x_{2})w\\
&=&(x_{0}+x_{2})^{l}Y_{W}(Y_{\sharp}(u\otimes g,x_{0})(v\otimes h),x_{2})w.
\end{eqnarray*}
This proves that $W$ is a $V\sharp H$-module.
\end{proof}

\section{Realizing vertex algebras $M(1)$ and $V_{L}$ through smash products}
In this section, we realize the vertex (operator) algebras 
associated with infinite-dimensional Heisenberg Lie algebras 
and the vertex algebras associated with nondegenerate even lattices 
through the smash product construction. 

First we consider the vertex operator algebras 
associated with Heisenberg Lie algebras.
Let ${\bf h}$ be a finite-dimensional vector space equipped with 
a nondegenerate symmetric bilinear form $\<\cdot,\cdot\>$.
Viewing ${\bf h}$ as an abelian Lie algebra we have
the affine Lie algebra 
$\hat{\bf h}={\bf h} \otimes \C[t,t^{-1}]\oplus \C c$, where
$$[a\otimes t^{m},b\otimes t^{n}]=m\delta_{m+n,0}\<a,b\>c$$
for $a,b\in {\bf h},\; m,n\in \Z$, and $c$ is a central element.
For $h\in {\bf h},\; n\in \Z$, 
following the tradition we also use $h(n)$ for $h\otimes t^{n}$. 
Set
$$\hat{\bf h}^{-}={\bf h}\otimes t^{-1}\C[t^{-1}],$$
an abelian subalgebra of $\hat{\bf h}$.
Let $L(-1)$ be the linear endomorphism of 
$\hat{\bf h}^{-}$ defined by
\begin{eqnarray}
L(-1)(h(-n))=nh(-n-1)\;\;\;\mbox{ for }h\in {\bf h},\; n\ge 1.
\end{eqnarray}
Of course, $L(-1)$ is a derivation of $\hat{\bf h}^{-}$ as a Lie algebra.
Set
\begin{eqnarray}
B_{\bf h}=U(\hat{\bf h}^{-})=S(\hat{\bf h}^{-}),
\end{eqnarray}
which is a vertex bialgebra by Example \ref{exahopf-Lie}.

Next we are going to define a $B_{\bf h}$-module vertex algebra structure on
$B_{\bf h}$. Notice that the adjoint module structure of $B_{\bf h}$ does not make
$B_{\bf h}$ a $B_{\bf h}$-module vertex algebra.

For $\alpha\in {\bf h}$, set $\alpha(0)=0$ on $B_{\bf h}$.
For $n\ge 1$, define $\alpha(n)$ to be the derivation of $B_{\bf h}$ with
\begin{eqnarray}
\alpha(n)(\beta(-m))=n\<\alpha,\beta\>\delta_{n,m}
\;\;\;\mbox{ for }\beta\in {\bf h},\; m\ge 1.
\end{eqnarray}
Set
\begin{eqnarray}
\alpha(x)^{-}=\sum_{n\ge 0}\alpha(n)x^{-n-1}\in (\Der B_{\bf h})[[x^{-1}]].
\end{eqnarray}

\bl{l-heisenberg-prepare}
For $\alpha\in {\bf h}$, we have 
$\alpha(x)^{-}\in {\rm PDer}^{-}(B_{\bf h})$.
There exists a $B_{\bf h}$-module structure
$Y_{M}$ on $B_{\bf h}$, uniquely determined by
$$Y_{M}(\alpha(-1),x)=\alpha(x)^{-}\;\;\;\mbox{  for }\alpha\in {\bf h}.$$
Furthermore, $B_{\bf h}$ equipped with this $B_{\bf h}$-module structure
$Y_{M}$ is a $B_{\bf h}$-module vertex algebra.
\el

\begin{proof} For $\alpha\in {\bf h},\; u\in B_{\bf h}$, since $\alpha(n)u=0$ 
for $n$ sufficiently large, we have 
$$\alpha(x)^{-}u\in B_{\bf h}[x^{-1}]\subset B_{\bf h}\otimes \C((x)).$$
For $\alpha\in {\bf h},\; n\ge 1$, with $L(-1)$ and $\alpha(n)$ being derivations of
$B_{\bf h}$, it is straightforward to show that
$$[L(-1),\alpha(n)]=-n\alpha(n-1).$$
That is, $[L(-1),\alpha(x)^{-}]=\frac{d}{dx}\alpha(x)^{-}$.
By Proposition \ref{pclassical-3} we have
$\alpha(x)^{-}\in {\rm PDer}^{-}(B_{\bf h})$.
For $\alpha,\beta\in L,\; h\in {\bf h},\; m,n\ge 0$, it is also 
straightforward to show that
$$\alpha(n)\beta(m)=\beta(m)\alpha(n).$$
Then $\alpha(x)^{-}$ for $\alpha\in {\bf h}$ generate a commutative differential algebra
in $B(B_{\bf h})$.  
As $B_{\bf h}=S(\hat{\bf h}^{-})$ is free, there exists a unique
differential algebra homomorphism from $B_{\bf h}$ to $B(B_{\bf h})$,
sending $\alpha(-1)$ to $\alpha(x)^{-}$ for $\alpha\in {\bf h}$.  It
follows from Lemma \ref{lmodule-algebra} that this homomorphism makes
$B_{\bf h}$ a $B_{\bf h}$-module vertex algebra with
$Y_{M}(\alpha(-1),x)=\alpha(x)^{-}$ for $\alpha\in {\bf h}$.  Since
$B_{\bf h}$ as a vertex algebra is generated by $\alpha(-1)$ for
$\alpha\in {\bf h}$, the uniqueness is clear.
\end{proof}

Now we fix this $B_{\bf h}$-module vertex algebra structure.
Then we have the smash product nonlocal vertex algebra 
$B_{\bf h}\sharp B_{\bf h}$. For $\alpha\in {\bf h}$, set
\begin{eqnarray}
\alpha(x)^{+}=\sum_{n\ge 1}\alpha(-n)x^{n-1}\in S(\hat{\bf h}^{-})[[x]]=B_{\bf h}[[x]].
\end{eqnarray}
Recall that there exists a unique vertex algebra structure 
$Y_{M(1)}$ on $S(\hat{\bf h}^{-})$ with 
\begin{eqnarray}
Y_{M(1)}(\alpha(-1),x)=\alpha(x)^{+}+\alpha(x)^{-}=\alpha(x)
\;\;\;\mbox{  for }\alpha\in {\bf h}.
\end{eqnarray}
This vertex algebra is denoted by $M_{\hat{\bf h}}(1,0)$ or by $M(1)$
for short.

\bp{pdouble-heisenberg}
The associative subalgebra $\Delta(B_{\bf h})$ of 
$B_{\bf h}\otimes B_{\bf h}$ 
is a vertex subalgebra of $B_{\bf h}\sharp B_{\bf h}$. 
Furthermore, $\Delta$ viewed as a map from
$M_{\hat{\bf h}}(1,0)$ to $\Delta(B_{\bf h})$,
\begin{eqnarray}
\Delta: M_{\hat{\bf h}}(1,0)=S(\hat{\bf h}^{-})=B_{\bf h}
\rightarrow \Delta(B_{\bf h})\subset 
B_{\bf h}\otimes B_{\bf h}=B_{\bf h}\sharp B_{\bf h},
\end{eqnarray}
is a vertex algebra isomorphism.
\ep

\begin{proof} For $\alpha\in {\bf h}$, we have
\begin{eqnarray}
Y_{B_{\bf h}}(\alpha(-1),x)=e^{xL(-1)}\alpha(-1)
=\sum_{n\ge 0}\alpha(-n-1)x^{n}=\alpha(x)^{+}.
\end{eqnarray}
Then
\begin{eqnarray}\label{edelta-alpha}
Y^{\sharp}(\Delta(\alpha(-1)),x)
&=&Y^{\sharp}(\alpha(-1)\otimes 1+1\otimes \alpha(-1),x)\nonumber\\
&=&Y_{B_{\bf h}}(\alpha(-1),x)\otimes 1+\alpha(x)^{-}\otimes 1
+1\otimes Y_{B_{\bf h}}(\alpha(-1),x)\nonumber\\
&=&\alpha(x)^{+}\otimes 1+\alpha(x)^{-}\otimes 1+1\otimes \alpha(x)^{+}.
\end{eqnarray}
Writing $Y^{\sharp}(u,x)=\sum_{n\in \Z}u_{n}^{\sharp}$ for $u\in B_{\bf h}\sharp B_{\bf h}$, 
we have
\begin{eqnarray}
& &\Delta(\alpha(-1))^{\sharp}_{-n}=\alpha(-n)\otimes 1+1\otimes \alpha(-n)=\Delta(\alpha(-n))
\;\;\;\mbox{ for }n\ge 1,\label{edelta-negative}\\
& &\Delta(\alpha(-1))^{\sharp}_{m}=\alpha(m)\otimes 1\;\;\;\mbox{ for }m\ge 0.
\label{edelta-positive}
\end{eqnarray}
Furthermore, for $\alpha,\beta\in {\bf h},\; m,n\in \Z$, we have
\begin{eqnarray}\label{esharp-heisenberg}
[\Delta(\alpha(-1))^{\sharp}_{m},\Delta(\beta(-1))^{\sharp}_{n}]
=m\<\alpha,\beta\>\delta_{m+n,0}.
\end{eqnarray}
Then $B_{\bf h}\sharp B_{\bf h}$ is an $\hat{\bf h}$-module of level $1$ with
$\alpha(m)$ acting as $\Delta(\alpha(-1))^{\sharp}_{m}$ 
for $\alpha\in {\bf h},\; m\in \Z$.
Since $B_{\bf h}$ as an algebra is generated by $\alpha(-n)$ for $\alpha\in {\bf h},\; n\ge 1$
and since $\Delta$ is an algebra homomorphism, 
it follows the P-B-W theorem and from (\ref{edelta-negative}) that 
$\Delta(B_{\bf h})$ is exactly the $\hat{\bf h}$-submodule of
$B_{\bf h}\sharp B_{\bf h}$ generated by ${\bf 1}\otimes {\bf 1}$. 
It then follows that $\Delta(B_{\bf h})$ is exactly the nonlocal vertex subalgebra
of $B_{\bf h}\sharp B_{\bf h}$ generated by $\Delta(\alpha(-1))$ for $\alpha\in {\bf h}$.
In view of the P-B-W theorem, $\Delta$ is injective, so that
$\Delta$, viewed as a map from $M_{\hat{\bf h}}(1,0)$ to
$\Delta(B_{\bf h})$, is an isomorphism of $\hat{\bf h}^{-}$-modules. 
Then from (\ref{esharp-heisenberg}), $\Delta$, viewed as a map from $M_{\hat{\bf h}}(1,0)$ to
$\Delta(B_{\bf h})$, is an isomorphism of $\hat{\bf h}$-modules. 
Consequently, $\Delta$ is an isomorphism of vertex algebras.
\end{proof}

We next consider the vertex algebras associated with nondegenerate even lattices.
Let $L$ be a nondegenerate even lattice;
a free abelian group of finite rank equipped with a nondegenerate $\Z$-valued 
symmetric $\Z$-bilinear form
$\<\cdot,\cdot\>$ such that $\<\alpha,\alpha\>\in 2\Z$ for $\alpha\in L$. 
Set ${\bf h}=\C\otimes_{\Z}L$ and then linearly
extend the form on $L$ to ${\bf h}$. As in the previous example we have
an affine Lie algebra $\hat{\bf h}$ and a differential bialgebra $B_{\bf h}$.
Let $\epsilon: L\times L\rightarrow \C^{\times}$ be a map such that
\begin{eqnarray*}
& &\epsilon(\alpha,0)=\epsilon(0,\alpha)=1,\\
& &\epsilon(\alpha,\beta+\gamma)\epsilon(\beta,\gamma)
=\epsilon(\alpha+\beta,\gamma)\varepsilon(\alpha,\beta)
\end{eqnarray*}
for $\alpha,\beta,\gamma\in L$. Such a map is called 
a normalized $\C^{\times}$-valued $2$-cocycle of $L$ (as an abelian group).
Denote by $\C_{\epsilon}[L]$ the $\epsilon$-twisted group algebra of
$L$ which has a $\C$-basis indexed by $L$, say $\{ e_{\alpha}\;|\; \alpha\in L\}$,
with multiplication
$$e_{\alpha}e_{\beta}=\epsilon(\alpha,\beta)e_{\alpha+\beta}
\;\;\;\mbox{ for }\alpha,\beta\in L.$$ 
Set
\begin{eqnarray}
B_{L,\epsilon}=\C_{\epsilon}[L]\otimes B_{\bf h},
\end{eqnarray}
an associative algebra. 
Extend the derivation $L(-1)$ of $B_{\bf h}$ to a linear endomorphism of $B_{L,\epsilon}$ by 
\begin{eqnarray}
L(-1)(e_{\alpha}\otimes u)
=e_{\alpha}\otimes \alpha(-1)u+e_{\alpha}\otimes L(-1)u
\end{eqnarray}
for $\alpha\in L,\; u\in S(\hat{\bf h}^{-})$.  It is straightforward
to check that $L(-1)$ is a derivation of $B_{L,\epsilon}$.  This makes
$B_{L,\epsilon}$ a differential algebra and hence a nonlocal vertex
algebra.

Recall that we have operators $h(n)$ on $B_{\bf h}$ 
for $h\in {\bf h},\;n\in \Z$. We extend the actions of $h(n)$ on $B_{\bf h}$ 
to $B_{L,\epsilon}$ by
\begin{eqnarray}
& &h(n)(e_{\alpha}\otimes u)=e_{\alpha}\otimes h(n)u
\;\;\;\mbox{ for }n\ne 0,\\
& &h(0) (e_{\alpha}\otimes u)=\<h,\alpha\> e_{\alpha}\otimes u
\end{eqnarray}
for $\alpha\in L,\; u\in B_{\bf h}$.

For $\alpha\in {\bf h}$, following \cite{flm} set
\begin{eqnarray}
E^{\pm}(\alpha,x)
=\exp\left(\sum_{n\in \pm \Z_{+}}\frac{\alpha(n)}{n}x^{-n}\right).
\end{eqnarray}
We have
\begin{eqnarray}
& &[L(-1),E^{\pm}(\alpha,x)]
=-\left(\sum_{n\in \pm \Z_{+}}\alpha(n-1)x^{-n}\right)E^{\pm}(\alpha,x),
\label{eproperty1}\\
& &\frac{d}{dx}E^{\pm}(\alpha,x)
=-\left(\sum_{n\in \pm \Z_{+}}\alpha(n)x^{-n-1}\right)E^{\pm}(\alpha,x).
\label{eproperty2}
\end{eqnarray}

\bl{llattice-half}
For $\alpha\in L$, we have
\begin{eqnarray}
Y_{B_{L,\epsilon}}(e_{\alpha},x)=E^{-}(-\alpha,x)e_{\alpha},
\end{eqnarray}
where $Y_{B_{L,\epsilon}}$ denotes the vertex operator map of
the nonlocal vertex algebra $B_{L,\epsilon}$. 
\el

\begin{proof}
For $\alpha\in L$, as
$Y_{B_{L,\epsilon}}(e_{\alpha},x)=e^{xL(-1)}e_{\alpha}$,
we must prove $e^{-xL(-1)}E^{-}(-\alpha,x)e_{\alpha}=e_{\alpha}$.
Clearly, it is true for $x=0$. 
Using the properties (\ref{eproperty1}) and (\ref{eproperty2}) we have
\begin{eqnarray*}
& &\frac{d}{dx}\left(e^{-xL(-1)}E^{-}(-\alpha,x)e_{\alpha}\right)\\
&=&-e^{-xL(-1)}L(-1)E^{-}(-\alpha,x)e_{\alpha}
+e^{-xL(-1)}\left(\sum_{n\in -\Z_{+}}\alpha(n)x^{-n-1}\right)E^{-}(-\alpha,x)e_{\alpha}\\
&=&-e^{-xL(-1)}E^{-}(-\alpha,x)L(-1)e_{\alpha}+e^{-xL(-1)}E^{-}(-\alpha,x)\alpha(-1)e_{\alpha}\\
&=&0,
\end{eqnarray*}
as $L(-1)e_{\alpha}=\alpha(-1)e_{\alpha}$. Thus 
$e^{-xL(-1)}E^{-}(-\alpha,x)e_{\alpha}=e_{\alpha}$.
\end{proof}

Set 
\begin{eqnarray}
B_{L}=\C[L]\otimes B_{\bf h}=\C[L]\otimes S(\hat{\bf h}^{-}),
\end{eqnarray}
a Hopf algebra. Note that $B_{L}=B_{L,\epsilon}$ with $\epsilon$ being trivial.
It is straightforward to check that
the comultiplication $\Delta$ and the counit $\varepsilon$ are
homomorphisms of differential algebras.
Then $B_{L}$ is naturally a differential bialgebra.
(Note that for a general $\epsilon$, $B_{L,\epsilon}$ is not a bialgebra 
with $\Delta$ defined in the obvious way.)

Next, we define a $B_{L}$-module nonlocal vertex algebra structure
on $B_{L,\epsilon}$. For this purpose we need the following universal property
of $B_{L}$:

\bl{lH(L)-universal} 
Let $L$ be a nondegenerate even lattice. 
Let $A$ be any commutative associative algebra with a derivation $\partial$ and let
$f: \C[L]\rightarrow A$ be any homomorphism of algebras. Then $f$ can
be extended uniquely to a homomorphism of differential algebras from
$B_{L}$ to $A$.  
\el

\begin{proof} For $\alpha,\beta\in L$, we have
\begin{eqnarray*}
f(e^{-\alpha-\beta})\partial f(e^{\alpha+\beta})
&=&f(e^{-\alpha-\beta})\partial (f(e^{\alpha})f(e^{\beta}))\\
&=&f(e^{-\alpha-\beta})(f(e^{\beta})\partial f(e^{\alpha})
+f(e^{\alpha})\partial f(e^{\beta}))\\
&=&f(e^{-\alpha})\partial f(e^{\alpha})
+f(e^{-\beta})\partial f(e^{\beta}).
\end{eqnarray*}
In view of this, we have a linear map from ${\bf h}$ to $A$, 
sending $\alpha\in L$ to 
$f(e^{-\alpha})\partial f(e^{\alpha})$.
Furthermore, we have an algebra homomorphism $g$ from 
$S(\hat{\bf h}^{-})$ to $A$ such that 
$$g(\alpha(-1-n))=\frac{1}{n!}\partial^{n}
(f(e^{-\alpha})\partial f(e^{\alpha}))
\;\;\;\mbox{ for }\alpha\in L,\; n\in \N.$$
Then $f\otimes g$ is an algebra homomorphism from $B_{L}$ to $A$,
extending $f$. Now we show that 
$(f\otimes g)L(-1)=\partial (f\otimes g)$, that is, $f\otimes g$ is 
a homomorphism of differential algebras.
For $\alpha\in L,\; n\in \N$, 
since $\alpha(-n-1)=\frac{1}{n!}L(-1)^{n}\alpha(-1)$,
we have $g(L(-1)^{n}\alpha(-1))=\partial^{n}g(\alpha(-1))$.
Assume that $a,b\in B_{\bf h}$ such that
$gL(-1)(a)=\partial g(a)$ and $gL(-1)(b)=\partial g(b)$. Then
$$gL(-1)(ab)=g(aL(-1)b)+g(bL(-1)a)=g(a)\partial g(b)+g(b)\partial (a)
=\partial (g(a)g(b))=\partial g(ab).$$
It follows from induction that $gL(-1)=\partial g$.
Furthermore, for $\alpha\in L$, we have
$$(f\otimes g)L(-1)e^{\alpha}=(f\otimes g)(e^{\alpha}\otimes \alpha(-1))
=f(e^{\alpha})f(e^{-\alpha})\partial f(e^{\alpha})=\partial f(e^{\alpha})
=\partial(f\otimes g)e^{\alpha}.$$
Then it follows that $(f\otimes g)L(-1)=\partial (f\otimes g)$.
As $\C[L]$ generates $B_{\bf h}$ as a differential algebra,
the uniqueness is clear.
\end{proof}

\bp{pprepare-lattice-voa}
We have
$$E^{+}(-\alpha,x)x^{\alpha(0)}\in {\rm PEnd}^{-}(B_{L,\epsilon})\;\;\;\mbox{ for }\alpha\in L.$$
Furthermore, there exists a unique $B_{L}$-module structure $Y_{M}$ on $B_{L,\epsilon}$ 
such that
\begin{eqnarray}
Y_{M}(e^{\alpha},x)=E^{+}(-\alpha,x)x^{\alpha(0)}
\;\;\;\mbox{ for }\alpha\in L
\end{eqnarray}
and $B_{L,\epsilon}$ equipped with this $B_{L}$-module structure $Y_{M}$
is a $B_{L}$-module nonlocal vertex algebra.
\ep

\begin{proof} Set $\Phi_{\alpha}(x)=E^{+}(-\alpha,x)x^{\alpha(0)}$ for $\alpha \in L$.
It is straightforward to check that
for $\alpha\in L,\; n\ge 0$, $\alpha(n)\in \Der(B_{L,\epsilon})$ and 
$[L(-1),\alpha(n)]=-n\alpha(n-1)$.
Then 
$$[L(-1),\log \Phi_{\alpha}(x)]=\frac{d}{dx}\log \Phi_{\alpha}(x).$$
In view of Proposition \ref{pclassical-3}, 
we have 
$$\log \Phi_{\alpha}(x)\in {\rm Der}(B_{L,\epsilon}, 
B_{L,\epsilon}\otimes (\C((x))[\log x],-d/dx)).$$
Consequently, $\Phi_{\alpha}(x)\in {\rm PEnd}^{-}(B_{L,\epsilon})$.
Clearly, we have $\Phi_{0}(x)=1$ and
$$\Phi_{\alpha}(x)\Phi_{\beta}(x) =\Phi_{\alpha+\beta}(x)\;\;\;\mbox{
for }\alpha,\beta\in L.$$ 
Denote by $A$ the subalgebra of $B(B_{L,\epsilon})$ generated by
$(d/dx)^{n}\Phi_{\alpha}(x)$ for $n\in \N,\; \alpha\in L$. It is clear
that $A$ is a commutative differential subalgebra.
By Lemma \ref{lH(L)-universal}, there exists a 
homomorphism $\pi$ of differential algebras from $B_{L}$ to $A$ such
that $\pi(e^{\alpha})=\Phi_{\alpha}(x)$ for $\alpha\in L$. 
As $B_{L,\varepsilon}$ is a module for $B(B_{L,\varepsilon})$ as a nonlocal vertex algebra
(by Proposition \ref{pbv-algebra}), the homomorphism $\pi$ 
gives rise to a $B_{L}$-module structure $Y_{M}$ on
$B_{L,\varepsilon}$ such that $Y_{M}(e^{\alpha},x)=\Phi_{\alpha}(x)$ for $\alpha\in L$.
It follows from Lemma \ref{lmodule-algebra} that
$B_{L,\epsilon}$ equipped with the $B_{L}$-module structure $Y_{M}$ is 
a $B_{L}$-module nonlocal vertex algebra.
As $\C[L]$ generates $B_{L}$ as a vertex algebra, the uniqueness is clear.
\end{proof}

Now we fix the $B_{L}$-module nonlocal vertex algebra structure on $B_{L,\epsilon}$.
We have the smash product nonlocal vertex algebra $B_{L,\epsilon}\sharp B_{L}$.

\br{rvoa-lattices}
{\em Recall from \cite{flm} that for the vertex algebra $V_{L}$ 
associated with $L$, we have
$V_{L}=B_{L,\epsilon}=\C_{\epsilon}[L]\otimes S(\hat{\bf h}^{-})$ 
as a vector space, where $\epsilon$ is a normalized $2$-cocycle of $L$
satisfying the condition
\begin{eqnarray}
\epsilon (\alpha,\beta)\epsilon (\beta,\alpha)^{-1}=(-1)^{\<\alpha,\beta\>}
\;\;\;\mbox{ for }\alpha,\beta\in L.
\end{eqnarray}
Denote the vertex operator map of $V_{L}$ by $Y_{V_{L}}$. For $\alpha\in L$,
we have
\begin{eqnarray}
Y_{V_{L}}(e_{\alpha},x)=E^{-}(-\alpha,x)E^{+}(-\alpha,x)e_{\alpha}.
x^{\alpha(0)}
\end{eqnarray}}
%For $\alpha,\beta\in L$,  the following relation holds
%\begin{eqnarray}\label{elocality-ea-eb}
%(x_{1}-x_{2})^{\<\alpha,\beta\>}Y_{V_{L}}(e_{\alpha},x_{1})
%Y_{V_{L}}(e_{\beta},x_{2})
%=(-x_{2}+x_{1})^{\<\alpha,\beta\>}Y_{V_{L}}(e_{\beta},x_{2})
%Y_{V_{L}}(e_{\alpha},x_{1}).
%\end{eqnarray}}
\er

\bp{plattice} Let $L$ be a nondegenerate even lattice
and let $\epsilon$ be the normalized $2$-cocycle of $L$, 
which was used in the construction of the vertex algebra $V_{L}$. 
Set
\begin{eqnarray}
U=\coprod_{\alpha\in L}
\C (e_{\alpha}\otimes e^{\alpha})\otimes \Delta(B_{\bf h})
\subset B_{L,\epsilon}\sharp B_{L}.
\end{eqnarray}
Then $U$ is an ordinary vertex
subalgebra of $B_{L,\epsilon}\; \sharp \; B_{L}$ and the linear map
\begin{eqnarray}
\pi: V_{L}\rightarrow U;\ \ 
e_{\alpha}\otimes u\mapsto (e_{\alpha}\otimes
e^{\alpha})\otimes \Delta(u)
\end{eqnarray}
for $\alpha\in L,\; u\in S(\hat{\bf h}^{-})$
is an isomorphism of vertex algebras.
\ep

\begin{proof} With $\Delta$ being an isomorphism from
$B_{\bf h}$ to $B_{\bf h}\otimes B_{\bf h}$,
it is clear that $\pi$ is a linear isomorphism.
Then it remains to prove that $\pi$ is a homomorphism of vertex
algebras.
Let $\alpha,\beta\in L,\; u\in S(\hat{\bf h}^{-})=B_{\bf h}$.
{}From Remark \ref{rvoa-lattices} we have
\begin{eqnarray*}
Y_{V_{L}}(e_{\alpha},x)(e^{\beta}\otimes u)
=x^{\<\alpha,\beta\>}\epsilon(\alpha,\beta)(e_{\alpha+\beta}\otimes E^{-}(-\alpha,x)E^{+}(-\alpha,x)u).
\end{eqnarray*}
Using (\ref{edelta-negative}) and (\ref{edelta-positive}) we get
\begin{eqnarray}
& &\pi\left(Y_{V_{L}}(e_{\alpha},x)(e^{\beta}\otimes u)\right)\nonumber\\
&=&x^{\<\alpha,\beta\>}\epsilon(\alpha,\beta)(e_{\alpha+\beta}\otimes e^{\alpha+\beta})
\Delta\left(E^{-}(-\alpha,x)E^{+}(-\alpha,x)u\right)\nonumber\\
&=&x^{\<\alpha,\beta\>}\epsilon(\alpha,\beta)(e_{\alpha+\beta}\otimes e^{\alpha+\beta})
(E^{-}(-\alpha,x)E^{+}(-\alpha,x)\otimes E^{-}(-\alpha,x))\Delta(u).
\end{eqnarray}
On the other hand, as $\Delta(e^{\alpha})=e^{\alpha}\otimes e^{\alpha}$, 
using Lemma \ref{llattice-half} we have
\begin{eqnarray*}
Y^{\sharp}(e_{\alpha}\otimes e^{\alpha},x)
&=&Y_{B_{L,\epsilon}}(e_{\alpha},x)\Phi_{\alpha}(x)\otimes Y_{B_{L}}(e^{\alpha},x)\nonumber\\
&=&E^{-}(-\alpha,x)e_{\alpha}E^{+}(-\alpha,x)x^{\alpha(0)}\otimes E^{-}(-\alpha,x)e^{\alpha}.
\end{eqnarray*}
Then
\begin{eqnarray}
& &Y^{\sharp}(e_{\alpha}\otimes e^{\alpha},x)\psi(e_{\beta}\otimes u)\nonumber\\
&=&Y^{\sharp}(e_{\alpha}\otimes e^{\alpha},x)(e_{\beta}\otimes e^{\beta})\Delta(u)\nonumber\\
&=&x^{\<\alpha,\beta\>}\epsilon(\alpha,\beta) (e_{\alpha+\beta}\otimes e^{\alpha+\beta})
(E^{-}(-\alpha,x)E^{+}(-\alpha,x)\otimes E^{-}(-\alpha,x))\Delta(u).
\end{eqnarray}
Consequently, we have
\begin{eqnarray}
\pi\left(Y_{V_{L}}(e_{\alpha},x)(e^{\beta}\otimes u)\right)
=Y^{\sharp}(e_{\alpha}\otimes e^{\alpha},x)\pi(e_{\beta}\otimes u)
\end{eqnarray}
for $\alpha,\beta\in L,\; u\in S(\hat{\bf h}^{-})=B_{\bf h}$.
As $\C_{\epsilon}[L]$ generates $V_{L}$ as a vertex algebra,
$\pi$ is a homomorphism of nonlocal vertex algebras. 
Since $V_{L}$ is a vertex algebra,
$\pi$ is an isomorphism of vertex algebras.
\end{proof}

\br{rlattice-module}
{\em Let $P$ be the dual lattice of $L$, i.e.,
$$P=\{ h\in {\bf h}\;|\; \<\alpha,h\>\in \Z\;\;\;\mbox{ for all }\alpha\in L\}.$$
Extend $\epsilon$ to a map from $L\times P\rightarrow \C^{\times}$ such that
\begin{eqnarray}
\epsilon (\alpha,\beta)\epsilon(\alpha+\beta,\gamma)=
\epsilon(\alpha,\beta+\gamma)\epsilon(\beta,\gamma)
\end{eqnarray}
for $\alpha,\beta\in L,\; \gamma\in P$. 
Define a $\C_{\epsilon}[L]$-module structure on $\C[P]$ by
$$e_{\alpha}\cdot e^{\gamma}=\epsilon(\alpha,\gamma)e^{\alpha+\gamma}
\;\;\;\mbox{ for }\alpha\in L,\; \gamma\in P.$$
Set 
$$V_{P}=\C [P]\otimes B_{\bf h}=\C [P]\otimes S(\hat{\bf h}^{-}).$$
The space $V_{P}$ is naturally a $B_{L,\epsilon}$-module. On the other hand,
the same argument in Proposition \ref{pprepare-lattice-voa} shows that
there exists a (unique) $B_{L}$-module structure $Y_{M}$ on $V_{P}$ 
with $Y_{M}(e^{\alpha},x)=E^{+}(-\alpha,x)x^{\alpha(0)}$ for $\alpha\in L$.
Using the commutation relation between $E^{+}(\alpha,x_{1})$ and $E^{-}(\beta,x_{2})$
(see \cite{flm}) one sees that Proposition \ref{pmodule-smash} applies to our situation, so
$V_{P}$ is a $B_{L,\epsilon}\sharp B_{L}$-module.
As $V_{L}$ is a vertex subalgebra of $B_{L,\epsilon}\sharp B_{L}$,
$V_{P}$ is a $V_{L}$-module. This gives a new proof for the existence of
a $V_{L}$-structure on $V_{P}$ (cf. \cite{flm}, \cite{ll}).}
\er

\end{document}